\documentclass[letter,11pt]{article}

\usepackage{lipsum}
\usepackage{amsfonts,amsmath,bm}
\usepackage{graphicx}
\usepackage{epstopdf}
\usepackage{algorithmic}
\usepackage{wrapfig}
\usepackage{hyperref}
\usepackage[left=1in,top=1in,right=1in,bottom=1in,letterpaper]{geometry}
\ifpdf
  \DeclareGraphicsExtensions{.eps,.pdf,.png,.jpg}
\else
  \DeclareGraphicsExtensions{.eps}
\fi

\usepackage{amsopn}

\DeclareMathOperator{\proj}{Proj} 

\def\M{\mathcal{M}}
\def\C{\mathcal{C}}
\def\N{\mathcal{N}}

\def\V{\mathcal{V}}

\def\R{\mathcal{R}}

\def\P{\mathcal{P}}

\def\I{\mathcal{I}}
\def\T{\mathcal{T}}

\def\RR{\mathbb{R}}
\def\dx{~\mathrm{d} x}
\def\ud{\,\mathrm{d}}

\def\Cov{\mathrm{Cov}}
\graphicspath{{Figures/}}

\title{Point cloud discretization of Fokker-Planck
  operators for committor functions}

\author{Rongjie Lai
\thanks{Department of Mathematics, Rensselaer
    Polytechnic Institute, Troy, NY ({\tt lair@rpi.edu},
    \url{http://homepages.rpi.edu/\string~lair}).}  
    \and 
Jianfeng Lu
\thanks{Department of Mathematics, Department of Physics, and
    Department of Chemistry, Duke University, Durham, NC
    ({\tt jianfeng@math.duke.edu},
    \url{http://www.math.duke.edu/\string~jianfeng/}).}  }

\date{}

\begin{document}

\maketitle

\begin{abstract}
  The committor functions provide useful information to the
  understanding of transitions of a stochastic system between disjoint
  regions in phase space.  In this work, we develop a point cloud
  discretization for Fokker-Planck operators to numerically calculate
  the committor function, with the assumption that the transition
  occurs on an intrinsically low-dimensional manifold in the ambient
  potentially high dimensional configurational space of the stochastic
  system. Numerical examples on model systems validate the effectiveness of the proposed method. 
\end{abstract}



\section{Introduction}

The understanding of transition of a stochastic system between
disjoint regions in the phase space has applications in the study of
chemical reactions and thermally activated process (see  \cite{MetzlerOshaninRedner:14, EVa:10} and references therein). For such physical
processes, the underlying dynamics can be often described by a
stochastic process such as the overdamped Langevin equation (also
known as the Brownian dynamics), which we will focus on in this paper
(while our methods can be generalized to other scenarios):
\begin{equation}\label{eq:overdamp}
  \ud X_t = - \nabla U(X_t) \ud t + \sqrt{2 \beta^{-1}} \ud W_t,
\end{equation} 
where $X_t \in \Omega \subset \RR^{d}$ denotes the current
configuration of the system at time $t$, $U: \Omega \to \RR$ is a
potential function, $\beta = 1/(k_BT)$ is the inverse temperature
($k_B$ being the Boltzmann constant and $T$ the absolute temperature),
and $W_t$ is the standard $d$-dimensional Wiener process. Here
$\Omega$ is the configurational space of the system (with some
suitable boundary condition if $\Omega$ is not the full space). We are
interested in understanding the transition of the dynamics
\eqref{eq:overdamp} from a subset $A \subset \Omega$ (represents for
example the reactant state) to a disjoint subset $B \subset \Omega$
(represents for example the product state of the system). Direct
simulation of such transitions might be challenging since most of
these systems exhibit time scale separations: It takes much longer
time for the system to go from $A$ to $B$ compared to the intrinsic
time scale of the dynamics (which limits the time step size of
numerical simulation). As a result, the study of such rare transitions
requires novel methodology development, which has been a very active
area in physical chemistry and applied mathematics.

The transition path theory, proposed by E and Vanden-Eijnden in
\cite{EVa:06} and further developed in \cite{MeScVa:06, MeScVa:09,
  LuNolen:15, CaVa:14}, is a framework to study the rare transitions
in the phase space. See also the review article \cite{EVa:10}.  In the
transition path theory, the central role is played by the
\emph{committor function} $q$ \cite{Hummer:04, EVa:06}, which is the
probability that a trajectory starting from $x \in \Omega$ first hits
$B$ rather than $A$. Denote $\tau_{\Sigma}$ the first hitting time of a
subset $\Sigma \subset \Omega$, the committor function is defined by
\begin{equation}
  q(x) = \mathbb{P}_x(\tau_B < \tau_A).
\end{equation}
Intuitively the committor function indicates the progression of a
transition: It takes the value $q = 0$ on $A$, $q = 1$ on $B$ and
increases going from $A$ to $B$. In chemical terms, the committor
function can be understood as a reaction coordinate of the transition (see \cite{Peters:16} for a recent review on reaction coordinates).

It follows that $q$ solves the following PDE on
$\Omega \backslash (A \cup B)$ with Dirichlet boundary conditions
given on $A$ and $B$ \cite{EVa:06, LuNolen:15}:  
\begin{equation}\label{eq:committor}
  \begin{cases}
    L q = 0, & \text{in } \Omega \backslash (A \cup B); \\
    q = 0, & \text{in } A; \\
    q = 1, & \text{in } B,     
  \end{cases}
\end{equation}
where $L$ is the infinitesimal generator of the process \eqref{eq:overdamp}, given by 
\begin{equation}
  L = - \beta^{-1} \Delta + \nabla U \cdot \nabla. 
\end{equation}
Note that some boundary
conditions are needed in the above equation if $\Omega$ is not a closed manifold, which we
will come back to later.

Once the committor function is provided, we can easily obtain
information on the reaction rate, density of transition paths, current
of transition paths \cite{EVa:06, EVa:10, LuNolen:15}, which help
understand the stochastic system.  Moreover, we can also write down
the SDE for the transition path between $A$ and $B$ that only depends
on the committor function \cite{LuNolen:15}, which can be used for
transition path sampling \cite{DellagoBolhuisGeissler:02,
  BolhuisChandlerDellagoGeissler:02}. 

Solving the PDE \eqref{eq:committor} for the committor function is
however non-trivial due to the curse of dimensionality. As a result,
in the framework of transition path theory, further assumptions are
usually made to approximate the committor function: It is assumed that
the transition from $A$ to $B$ is concentrated in quasi-one
dimensional ``reaction tube'', which is the working assumption of the
finite temperature string method \cite{ERenVa:05, VaVe:09}. Due to the
usefulness of the committor function in understanding the transition,
other approaches have been also developed, mainly in the chemical
physics literature, for example methods based on statistical analysis
of an ensemble of trajectories \cite{MaDinner:05, PetersTrout:06,
  LechnerRogal:10}, based on approximating the diffusion by a discrete
state space jump process (milestoning) \cite{MajekElber:10,
  KirmizialtinElber:11} (see also \cite{CaVa:14} on direct computation
of committor function on a discrete Markov jump process).  To the best
of our knowledge, none of the existing methods aims at approximating
the PDE~\eqref{eq:committor} directly.

The main contribution of this work is to propose a method to directly
solve the committor equation \eqref{eq:committor} based on point cloud
discretization, especially the technique of local mesh discretization
recently developed in \cite{LaiLiangZhao:13}. Assume a given point cloud 
samples the equilibrium distribution of the dynamics
\eqref{eq:overdamp}, the idea is to discretize the PDE
\eqref{eq:committor} on the point cloud to approximate the committor
function. Here, the working assumption is that while the
configurational space of the stochastic system is high dimensional,
the transition between the interested regions $A$ and $B$ lies in an
intrinsically low-dimensional manifold (for simplicity, we assume that
the intrinsic dimension does not change).  In particular, this
generalizes the ``reaction tube'' assumption of the finite temperature
string method to transition in higher than quasi-one dimension.

Our method is closely related in spirit to the method of diffusion map
\cite{CoifmanLafon:06}. In particular, in the work \cite{NaLaCoKe:06,
  CoKeLaMaNa:08} and subsequently \cite{MaggioniClementi:11a,
  MaggioniClementi:11b}, the diffusion map has been applied to obtain
an approximation of the infinitesimal generator to compute the first
few low-lying eigenfunctions of $L$, the Fokker-Planck operator. Those
eigenmodes are used to approximate the long time behavior and to
coarse-grain the dynamics \eqref{eq:overdamp}. As demonstrated in
\cite{MaggioniClementi:11a, MaggioniClementi:11b}, the assumption that
the stochastic dynamics can be approximated by low-dimensional
reaction modes is indeed valid for a variety of chemical systems. It
is in fact possible to use the idea of diffusion maps to approximate
the equation for the committor function \eqref{eq:committor} and we
will compare our method with the diffusion map based method. Let us
also remark that we do not focus on the low-lying eigenmodes of
Fokker-Planck equation as \cite{NaLaCoKe:06, MaggioniClementi:11a,
  MaggioniClementi:11b}, but rather the committor function, which
provides the information on the transition from $A$ to $B$ more
directly.  One potential advantage of focusing directly at the
transition region is that we do not require point cloud to well
represent the regions $A$ and $B$, which are typically of high
intrinsic dimension, as usually they correspond to regions of local
minima of the potential energy surfaces, and hence leads to challenges
for point cloud discretization.  We refer the readers to \cite{EVa:10}
for further comparison between the committor function and low-lying
eigenmodes. Besides diffusion map and the local mesh method, other
numerical techniques have been also developed for solving PDEs on
point clouds without a mesh and this has been a very active area in
recent years; a brief review of those can be found in
\ref{sec:localmesh}.

The rest of the paper is organized as follows. The point cloud discretization method for the committor function is described in detail in \ref{sec:localmesh}. The algorithm is validated and compared with other approaches through numerical examples in \ref{sec:experiments}. We conclude the paper in \ref{sec:conclusions}.

\section{Local mesh method for Fokker-Planck operators}\label{sec:localmesh}

Given a potential function $U(x)$ defined on $\Omega$, recall that we
aim at solving the following equation on $\Omega$ represented as point
clouds:
\begin{equation}
\begin{cases}
-\beta^{-1} \Delta q(x) + \nabla U(x) \cdot \nabla q(x) = 0, & x\in \Omega \backslash (A\cup B);  \\
q(x) = 0, &  x\in A; \\
q(x) = 1, & x\in B;  \\
\nabla q(x)\cdot  \vec{n}(x) = 0, & x\in \partial\Omega. 
\end{cases}
\label{eqn:FK}
\end{equation}
Note that compared to \eqref{eq:committor}, we have also specified the
Neumann boundary condition at the boundary of $\partial\Omega$, which
is the natural boundary condition from a variational point of view, as
discussed below.

Before we proceed to numerical algorithms, let us write down the weak formulation of \eqref{eqn:FK}. Given any test function $\eta(x)$, using the equation and integration by parts, we have:
\begin{equation}
\begin{aligned}
0 & = \int_{\Omega}  \left(\Delta q(x) -\beta \nabla U(x) \cdot \nabla q(x) \right ) \eta(x) e^{-\beta U(x)}   \dx \\
 &  =  \int_{\Omega} \nabla \left(\nabla q(x) \eta(x) e^{-\beta U(x)} \right) - \nabla q(x) \nabla \eta(x) e^{-\beta U(x)} \dx  \\
 & =  - \int_{\Omega}\nabla q(x) \nabla \eta(x) e^{-\beta U(x)} \dx, 
 \label{eqn:weakFK}
\end{aligned}
\end{equation}
where in the last equality, we have used the Neumann boundary
condition $\vec{n}(x) \cdot \nabla q(x) = 0$ on $\partial \Omega$ so
that the boundary contribution vanishes. Note that the Gibbs weight
$e^{-\beta U(x)}$ appeared in \eqref{eqn:weakFK} is the invariant measure $\rho(x) = Z^{-1}e^{-\beta U(x)}$ of the overdamped
equation \eqref{eq:overdamp} up to a normalization constant $Z = \int_{\Omega} e^{-\beta U(x)} \dx$. 

In the transition path theory~\cite{EVa:10}, after we compute $q(x)$, we can immediately obtain the transition rate between $A$ and $B$ provided by 
\begin{equation}\label{eq:nuR}
  \nu_R = k_B T \int_{\Omega} \lvert \nabla q(x) \rvert^2 \rho(x) \dx, 
\end{equation}
which we will later use to quantify the approximation quality of our numerical methods. The committor function also leads to other useful information to understand the dynamics, including the probability density of reactive trajectories 
\begin{equation}
  \rho_R(x) = q(x) (1 - q(x)) \rho(x) 
\end{equation}
and also the reactive current
\begin{equation}
  J_R(x) = k_B T \rho(x) \nabla q(x),
\end{equation}
which can be obtained through explicit formulas based on $q$. 


Let us come back to numerical solution to \eqref{eqn:FK}. Since
$\Omega$ is potentially high dimensional while the transition between the interested regions A and B lies in an intrinsically low-dimensional manifold, our idea is to solve the
equation and get the approximation of $q$ instead on a point cloud well sampling the transition path.
More precisely, assume that we are given a point cloud
$\P=\{\bm{p}_i\in\mathbb{R}^N \mid i=1,\dots,n\}$, for instance from
snapshots of a long trajectory of the SDE
\eqref{eq:overdamp},\footnote{The advantage of using a long ergodic
  trajectory of the SDE is to cover the important region of the
  transition; on the other hand, let us emphasize that we do not need
  to assume that the point cloud is distributed according to the
  invariant measure of the SDE \eqref{eq:overdamp}, in particular, the
  method works as long as the point cloud well cover the important
  part of the configurational space.} the goal here is to approximate the
committor function, given by \eqref{eq:committor} based on $\P$. In particular, we
are aiming at the value of $q$ on the point clouds, instead of on the
whole configurational space $\Omega$. We will also interpret the sets
$A$ and $B$ as the collection of points in the point cloud $\P$ that
lies in the two chosen sets in the whole configurational space. For
simplicity of notation, when there is no danger of confusion, we will
not distinguish in the sequel the sets $A$, $B$, $\Omega$ etc. with
their point cloud interpretations.

In the rest of this section, we first give a brief review of the existing diffusion map based method~\cite{CoifmanLafon:06,MaggioniClementi:11a} for approximating 
the Fokker-Planck equation on the given point clouds, which will be compared later with our method. After that, our method of solving the Fokker-Planck equation on point clouds will be presented based on the local mesh method discussed in~\cite{LaiLiangZhao:13}.

\subsection{Diffusion map discretization of Fokker-Planck operator}

In recent years, several numerical schemes have been proposed for
solving partial differential equations on point clouds without a
global mesh or grid in the ambient space. These methods are
particularly useful in high dimensions where a global triangulation or
grid is intractable. Among those methods, a popular class of methods
is kernel based, where the Laplace-Beltrami operator on point clouds
is approximated by heat diffusion in the ambient Euclidean
space~\cite{CoifmanLafon:06} or in the tangent
space~\cite{Belkin09clp} among nearby points. In other words, the
metric on the manifold is approximated by Euclidean metric
locally. The main advantage of such methods is their
simplicity and generality for approximating diffusion type operators. 
Among those, the methods based on the ideas of diffusion map~\cite{CoifmanLafon:06} have been rather popular. For discretizing the Fokker-Planck operator $ -\beta^{-1}\Delta + \nabla U\cdot \nabla$ in \eqref{eq:committor} on the point cloud, one uses 
 $L = D^{-1} K - I_n$ where $K$ is the normalized kernel function according to the Gibbs weight $\exp(-\beta U)$:
\begin{equation}\label{eqn:diffKernel}
K(\bm{p}_i, \bm{p}_j) =\exp(-\beta U(\bm{p}_i))  \exp \left(-\frac{\|\bm{p}_i - \bm{p}_j\|^2}{2\epsilon_i\epsilon_j} \right)\exp(-\beta U(\bm{p}_j))
\end{equation}
and $D$ being a diagonal matrix with 
\[
D_{ii} = \sum_{j} K(\bm{p}_i,\bm{p}_j), 
\] 
and $I_n$ being the identity matrix of size $n \times n$. Note that $\epsilon_i$ can be chosen either as a constant or adaptively depending on the local information of the data set using the multiscale SVD method proposed in \cite{little2009multiscale,MaggioniClementi:11a}. 
Therefore, if we denote $C = \Omega - A\cup B$ (recall that these sets
are interpreted as the subsets of points in $\P$), then the committor
function $q = [q_A,q_B, q_C]$ can be solved by the linear equation
\[
L(C,C)q_C = - L(C, B) q_B,
\]
where $q_B \equiv 1$ and $q_A \equiv 0$ according to the boundary
conditions of the committor function. We remark that $L(C,C)$ and $L(C,B)$ are denoted as the restriction of the matrix $L$ on the index set of $(C,C)$ and $(C,B)$, respectively.

The above kernel based methods however usually render low order
approximation. More recently, two new
methods for solving PDEs on point clouds are introduced
in~\cite{Liang:CVPR2012,Liang2013solving,LaiLiangZhao:13}, where the
differential operators are approximated systematically and
intrinsically at each point from local approximation of the manifold
and its metric through nearby neighbors.
These methods can achieve
high order accuracy and can be used to approximate general differential
operators (\textit{i.e.}, other than diffusion type operators) on point clouds sampling manifolds with arbitrary dimensions
and co-dimensions.  Moreover, the computational complexity depends
mainly on the intrinsic dimension of the manifold rather than the dimension
of the embedded space.

To solve \eqref{eqn:FK} on a point cloud, we will focus on
generalizing the local mesh method discussed in \cite{LaiLiangZhao:13}
to Fokker-Planck operators. The local mesh method is natural to
handle the Neumann boundary conditions, which is more advantageous in
the current scenario. The crucial part of solving the weak formula \eqref{eqn:weakFK} is to numerically approximate $\int_\Omega \nabla q \nabla \eta e^{-\beta U} \mathrm{d} x$ on the given point cloud. Note that the definition of differential operators only rely on local information of point clouds. Thus, our idea is composed of constructing local connectivity to approximate the gradient operator and estimating the stiffness matrix on point clouds. In the rest of this section, we will discuss technical details about local connectivity construction and the stiffness matrix approximation on point clouds. After that, a numerical method of solving the Fokker-Planck equation \eqref{eqn:FK} will be proposed. 


\subsection{Local connectivity construction for point clouds}
\label{subsec:localmesh}

\begin{wrapfigure}{r}{6cm}
\vspace{-0.3cm}
\centering
\includegraphics[width=1\linewidth]{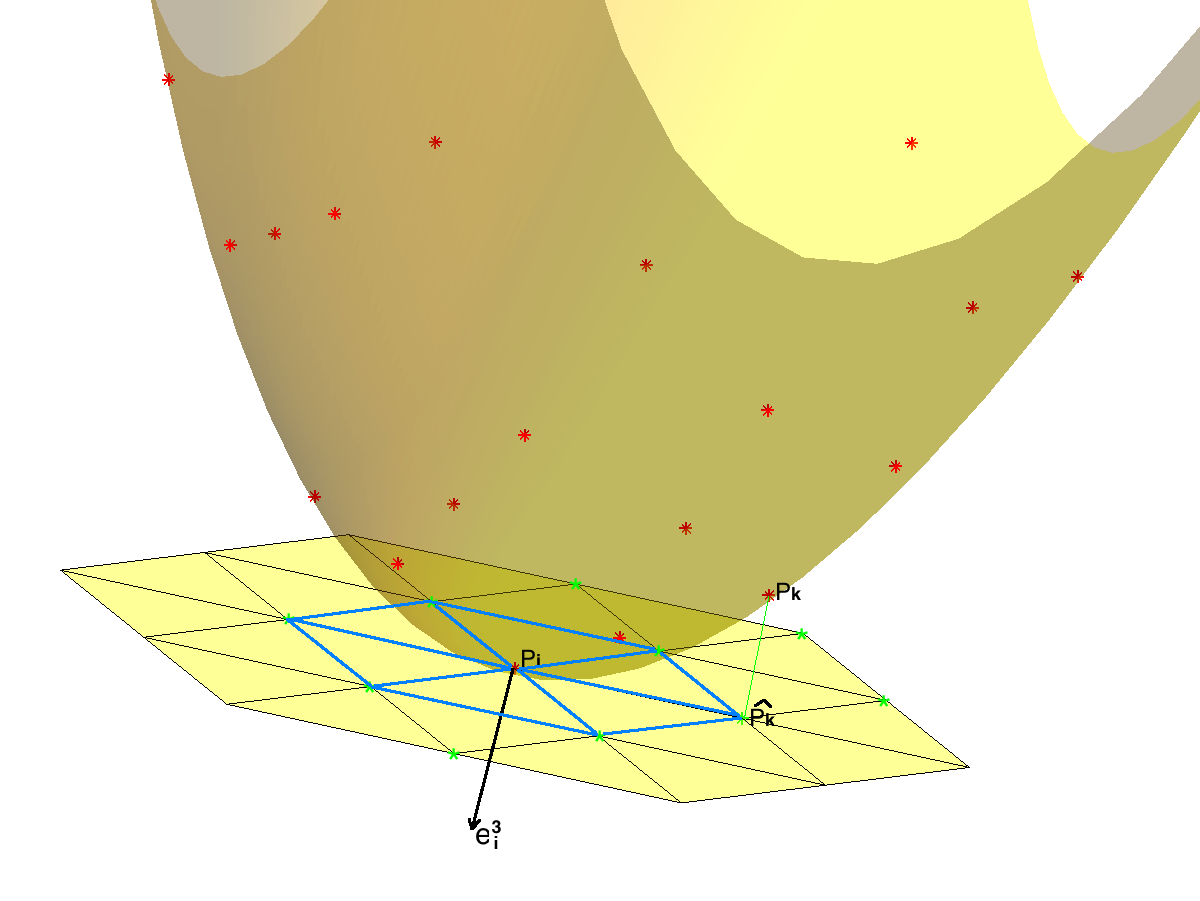}\\
\label{fig:PC_localmesh_pro}
\vspace{-.5cm}
\caption{Red stars mark the KNN of $\bm{p}_i$, green stars mark the projection of red stars on the tangent plane at $\bm{p}_i$, and blue triangles color-code the connectivity of the first ring of $\bm{p}_i$.}
\end{wrapfigure}
Given a point cloud $\P=\{\bm{p}_i\in\mathbb{R}^N~|~ i=1,\dots,n\}$ sampled from a $d$-dimensional manifold $\M$ in $\mathbb{R}^N$, let us denote the indices set of K-nearest neighborhood (KNN) of each
point $\bm{p}_i \in \P$ by $\I(i)$ and write $\P(i) = \{\bm{p}_k \in \P~|~ k\in \I(i)\}$. Based on $\P(i)$, the tangent space and normal space of $\M$ at $\bm{p}_i$ can be approximated by the standard principle component analysis (PCA)~\cite{jolliffe2002principal}. More precisely, one can first construct the following covariance matrix $\mathrm{Cov}_i$ of $\P(i)$:
\begin{equation}
  \mathrm{Cov}_i=\sum_{k\in \I(i)}(\bm{p}_k - \bm{c}_i)^T(\bm{p}_k-\bm{c}_i)   
\end{equation} 
where $\bm{c}_i$ is the local barycenter of $\P(i)$, given by 
\begin{equation*}
  \bm{c}_i=\frac{1}{|\P(i)|}\sum_{k\in \I(i)}\bm{p}_k. 
\end{equation*}
As the point cloud is sampled from a $d$-dimensional manifold $\M$, if the local sampling is dense enough to resolve local features, then eigenvalues of $\Cov_i$ have a natural jump $\lambda_i^1 \geq \cdots \geq \lambda_i^d \gg \lambda_i^{d+1} \geq \cdots\geq \lambda_i^N \geq 0$
which guides the splitting $\mathbb{R}^N = \T_i \oplus \N_i$. 
Here $\T_i$ represents the tangent space spanned by $\{\bm{e}_i^1,\cdots,\bm{e}_i^d\}$ corresponding to the $d$ largest eigenvectors of $\mathrm{Cov}_i$, and $\N_i$ represents the normal space spanned by $\{\bm{e}_i^{d+1},\cdots,\bm{e}_i^N\}$.
For noisy data set, a technique called multiscale SVD method~\cite{little2009multiscale} is an effective way to estimate the intrinsic dimension. 
To construct the local connectivity and a mesh at the point $\bm{p}_i$, we project its K-nearest neighborhood $\N(i)$ on the tangent plane $\T_i$ of $\bm{p}_i$. Namely, we have the following construction:
\begin{equation}
\hat{\bm{p}}_k = \proj_{\T_i}(\bm{p}_k) = \bm{p}_k - \bm{p}_i - \sum_{\alpha=d+1}^N\langle \bm{p}_k-\bm{p}_i, \bm{e}_i^\alpha\rangle \bm{e}_i^\alpha,\quad k\in\I(i)
\end{equation}
For convenience, we have
translated the K-nearest neighborhood $\P(i)$ to
center $\bm{p}_i$ at the origin. 
With this projection, all points in $\{\hat{\bm{p}}_k, k\in I(i)\}$ belong to the tangent plane $\T_i$. Then, the local mesh structure near $\bm{p}_i$ can be obtained by the standard Delaunay triangulation. Denote by $\R(i)=\{F_i^1,\cdots, F_i^{l_i}\}$ all simplexes adjacent to $\bm{p}_i$, referred as $\bm{p}_i$'s first ring and $\V(i)$ all vertices in the first ring of $\bm{p}_i$. The local connectivity of $\bm{p}_i$ is provided by $\C_i = \{\bm{p}_i; \V(i),\R(i)\}$. Figure~\ref{fig:PC_localmesh_pro} illustrates an example of a set of point sampled in a 2D manifold in $\mathbb{R}^3$.
After obtaining the local connectively, we are ready to discretize the weak formula \eqref{eqn:weakFK} to solve the Fokker-Planck equation \eqref{eqn:FK}. Our basic idea is to represent $q$ as a linear combination of a nodal basis on the point cloud. Derivatives of all nodal basis can be approximated by linear interpolation. After that, we approximate the stiffness matrix using the nodal basis and solve the weak equation. 

\subsection{Stiffness matrix construction on point clouds}
\label{sec:stiffM}
Suppose that we have a point cloud $\P$ with local connectivity $\C = \{\C_i = \{\bm{p}_i; \V(i),\R(i)\} ~|~ i = 1,\cdots,n\}$ constructed in section~\ref{subsec:localmesh}, we define the set of nodal basis $\{\eta_j\}_{j=1}^{n}$ on $\P$ as:
\begin{eqnarray}
\eta_j: \P \rightarrow \mathbb{R},\quad \eta_j(\bm{p}_i) =
\begin{cases} 1, & \text{if }\quad i=j \\ 
0, & \text{otherwise }
\end{cases}
\quad i,j = 1, \cdots, n
\label{eqn:linearelements}
\end{eqnarray}
Inspired by the idea from the finite element method, we consider $\{\eta_i\}_{i=1}^n$ as the set of linear elements defined on $\P$ given in (\ref{eqn:linearelements}), and write $E= \text{Span}_{\mathbb{R}}\{\eta_i\}_{i=1}^{n}$. We approximate the desired committor function $q= (q(\bm{p}_1),\cdots,q(\bm{p}_n))^T$ defined on $\P$ as $q=\sum_i^n v_i\eta_i$. 

Then the discrete version of the equation \eqref{eqn:weakFK} on $\P$
can be defined by the following weak formulation:
\begin{eqnarray}
\int_{\P}e^{-\beta U} \nabla_{\P}q \cdot \nabla_{\P}
 \eta  = 0,~\quad \forall~\eta\in E. 
 \label{formula:LBeigsdiscrete}
\end{eqnarray}  
Therefore, the essential part of
the above problem is a numerical approximation of the stiffness matrix $S = (S_{ij})_{n\times n}$ for the given point cloud $\P$ with:
\begin{equation}\label{eq:stiffness}
S_{ij} = \int_{\P}e^{-\beta U} \nabla_{\P}\eta_i\cdot\nabla_{\P}\eta_j.
\end{equation}

To clearly indicate the proposed method, we would like to assume the point cloud is sampled from a two dimensional manifold in $\mathbb{R}^N$.\footnote{We emphasize that the proposed idea can be easily adapted to other cases that the intrinsic dimension $\M$ is not two.} In this case, the local connectivity is a set of local triangle mesh structure.  Ideally, we would like $S$ to
be symmetric and non-negative definite, similar to the usual
properties of the stiffness matrix for a triangulated
surface~\cite{dziuk1988finite}. 
Unfortunately, these global properties
might not be possible to achieve due to the use of only local mesh in our discretization. To be more specific, the first ring structure of $\bm{p}_i$ is not
necessary compatible with the first ring structure of $\bm{p}_j$ although
$\bm{p}_j$ belongs to the first ring of $\bm{p}_i$.  To have a numerical
approximation of the stiffness matrix, we first define
\begin{eqnarray}
A_{ij}&=&\sum_{F\in\R(i)}\int_{F} e^{-\beta U}\nabla_{F}\eta_i\cdot\nabla_{F}\eta_j  \nonumber \\
&=& \sum_{F\in\R(i),[\bm{p}_i,\bm{p}_j]\in F} -\frac{1}{2}w_{ij}^F \cot\alpha^{F}_{ij}(p_i) , \quad i\neq j
\end{eqnarray}
where $\int_{F} e^{-\beta U}\nabla_{F}\eta_i\cdot\nabla_{F}\eta_j$ is computed by linear interpolation of $U, \eta_j$ and $\eta_j$ on $F$,$w_{ij}^F$ is the average of $e^{-\beta U}$ on $F$ and $\alpha^F_{ij}(\bm{p}_i)$ are the angles opposite to the edge connecting points $\bm{p}_i$ and $\bm{p}_j$ in the face $F$. 
Note that $A_{ij}$ may not be equal to $A_{ji}$ due to the possible incompatibility of the first ring structures of $\bm{p}_i$ and $\bm{p}_j$. One simple symmetrized definition of the stiffness matrix is the following:
\begin{equation}
S_{ij}=
\left\{
\begin{array}{ll}
\displaystyle\frac{1}{2}\left(A_{ij}+A_{ji}\right),  &\quad  \mbox{if}   \quad i\ne j  \vspace{0.2cm}\\
\displaystyle -\sum_{k\neq i} S_{ik}, &\quad  \mbox{if}   \quad  i = j
\end{array}\right.
\label{eqn:Stiffness1}
\end{equation}
The above definition of the diagonal elements is to enforce the
consistency condition, \textit{i.e.}, constant function is an
eigenfunction of $S$ with zero eigenvalue. In particular, if
all triangles in the first ring structure are acute, off-diagonal
elements are non-positive and diagonal elements
$S_{ii}=|\sum_{k\neq i} S_{ik}|$ are positive. Hence all eigenvalues
are real and non-negative. When the density of points is
reasonably uniform on $\M$, this definition of stiffness matrix works
quite well. However, when the density of points are non-uniformly as the data produced from the SDE \eqref{eq:overdamp} used in our experiments, the
first ring structure of $\bm{p}_i$ is more likely incompatible with the
first ring structure of neighboring points $\bm{p}_j$. To overcome this issue, we use a similar
strategy used in \cite{LaiLiangZhao:13} to construct $S$ as follows: 
\begin{equation}
S_{ij}=\left\{
\begin{array}{ll}
\max (A_{ij}, A_{ji}) & \mbox{if} \quad A_{ij}\le 0 \quad \mbox{and} \quad A_{ji}\le0 \\
\min(A_{ij}, A_{ji}) & \mbox{if} \quad A_{ij}\ge 0 \quad \mbox{and} \quad A_{ji}\ge0 \\
\min(A_{ij}, A_{ji}) & \mbox{if} \quad A_{ij}\cdot A_{ji} <0  \\
-\sum_{k\neq i} S_{ik} & \mbox{if}  \quad  i = j
\end{array}
\right.
\label{eqn:Stiffness2}
\end{equation}
As long as the stiffness matrix are constructed, we can approximate the committor function, the solution of the equation~\eqref{eqn:FK}, in the following way. 

Remember that we write $C = \Omega - A\cup B$, then we obtain the following discretization of the equation~\eqref{eqn:FK}. 
\begin{equation}
\left\{\begin{array}{c}
S(C, [A, B, C]) (q_A,~ q_B,~ q_C)^T = 0. \\ 
q_A = 0 \\
 q_B = 1
\end{array} \right. \label{eqn:FK_LM}
\end{equation}
This provides a matrix equation $$S(C,C)q_{C} =   -S_{C,B} q_B$$ which solves the committor function we desired.



\section{Numerical Experiments}
\label{sec:experiments}
In this section, we test the proposed method for committor functions
on several examples obtained from the stochastic differential equation
\eqref{eq:overdamp}. All experiments are impletmented by MATLAB in a
PC with a 32G RAM and a 2.7 GHz quad-core CPU.

\subsection{Comparison of methods for a 1D double well potential}
\label{subsec:diffusionmap}

We first conduct numerical experiments on 1D interval $[-1,~1]$ for
the standard double well potential $U = (x^2-1)^2$ with $\beta = 1$.
We choose the sets $A = [-1,~-1+0.1]$ and $B = [1-0.1, ~1]$.  We will
compare our method with the diffusion map discretization of the
Fokker-Planck operator, following the
works~\cite{CoifmanLafon:06,CoKeLaMaNa:08,MaggioniClementi:11a,MaggioniClementi:11b}. The
numerical experiments illustrate that the local mesh method achieves
better accuracy and robustness.

\begin{figure}[h]
\begin{minipage}{0.49\textwidth} 
\centering
\includegraphics[width=1\linewidth]{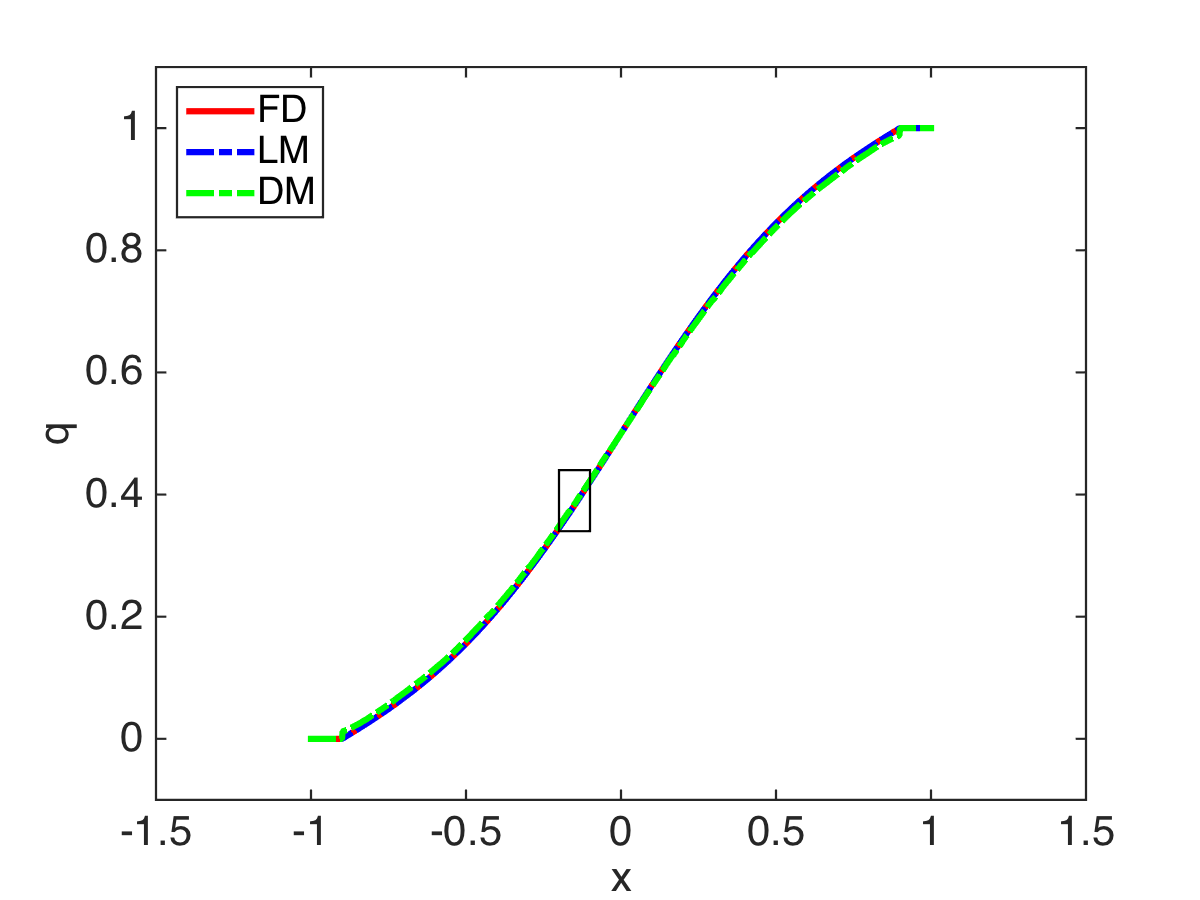}
\end{minipage}
\begin{minipage}{0.49\textwidth} 
\centering
\includegraphics[width=1\linewidth]{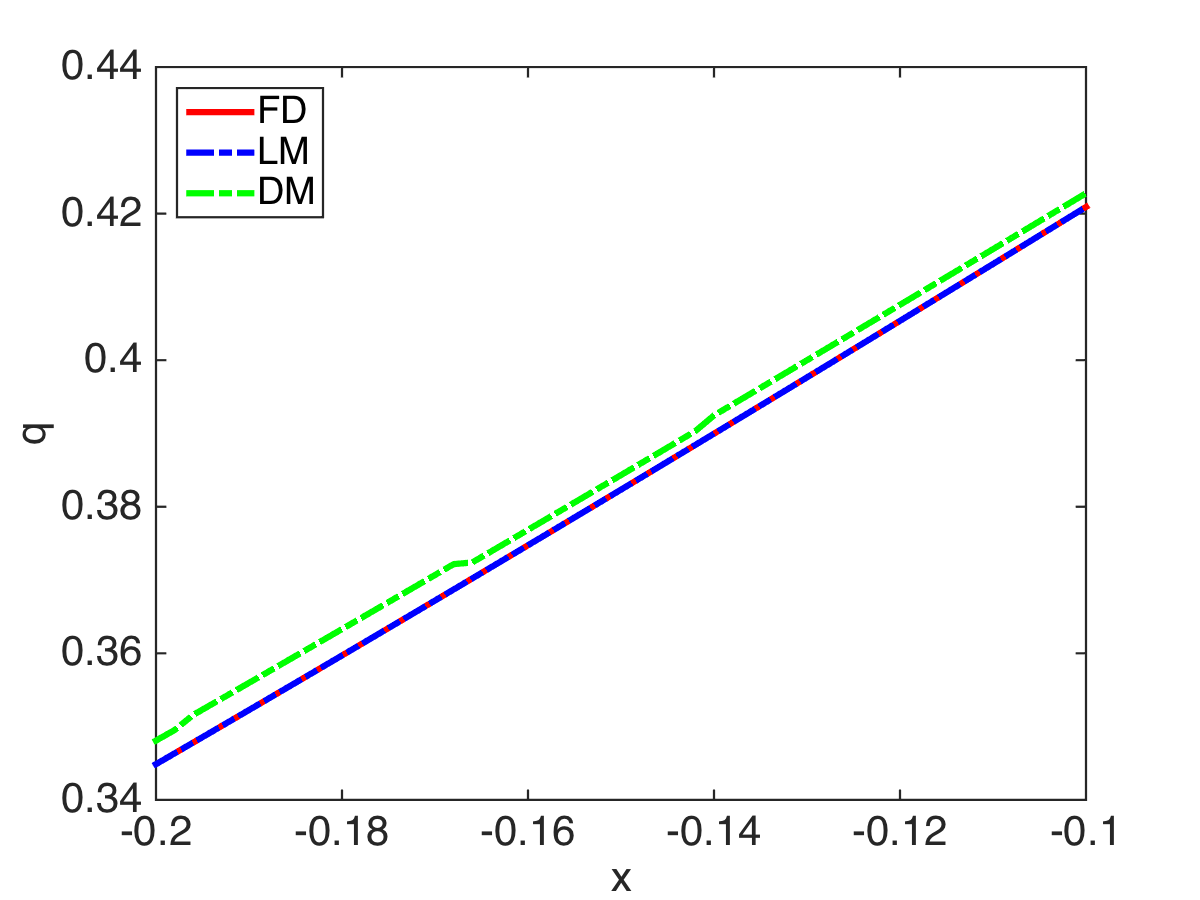}
\end{minipage}
\caption{Left: Solutions obtained from a finite element method (FE), the proposed local mesh method (LM) and the diffusion map method (DM). Right: A zoom-in image. }
\label{fig:Comparison_FD_LM_DM}
\end{figure}

We first compare results using our method and diffusion map method on
$1000$ equally sampled points distributed on $[-1,~1]$. In this case,
we also apply a standard finite element based method to solve the
weak equation \eqref{eqn:weakFK}, serving as a reference. Figure~\ref{fig:Comparison_FD_LM_DM}
illustrates comparisons among these three methods on the same point
cloud data set (for finite element, we interpret the points as grid
points). It is clear to see that solutions obtained from the finite
element method and the proposed method are nearly identical. The
diffusion map approach yields a similar result, however not quite as
accurate as the proposed method. To quantify this, we denote
$q_{\text{FE}}$, $q_{\text{LM}}$ and $q_{\text{DM}}$ as solutions
obtained from the finite element method, the proposed local mesh
method and the diffusion map method, respectively. The maximum
absolute error
\begin{equation}
  \max \left\{\frac{|q_{\text{FE}}(x) - q_{\text{LM}}(x)|}{|q_{\text{FE}}(x)|} \right\}
\end{equation}
between solutions using finite element and the proposed method is
$3.3317e$-$7$, while the maximum absolute error between solutions using
finite element and the diffusion map method is $0.0112$.

\begin{figure}[h]
\begin{minipage}{0.49\textwidth}
\centering
\includegraphics[width=1\linewidth]{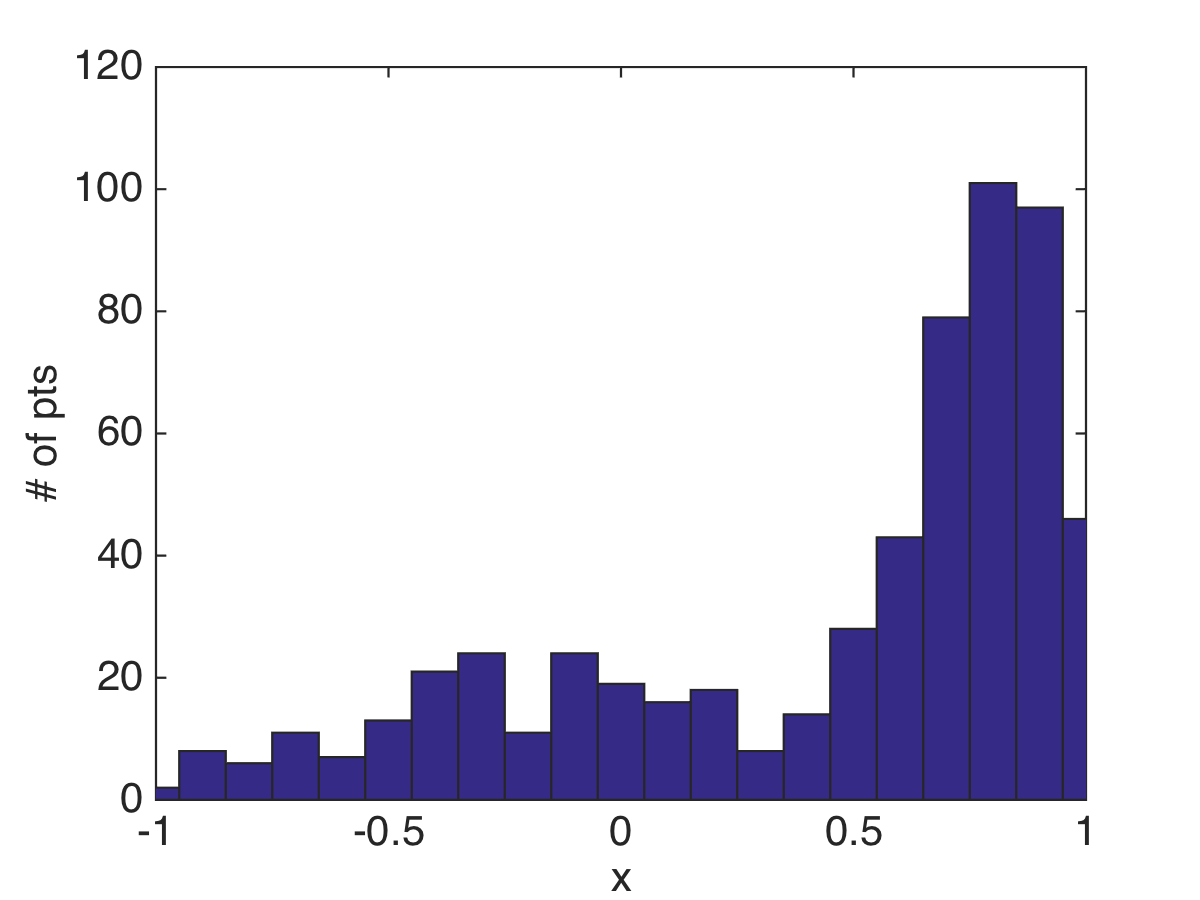}
\end{minipage}
\begin{minipage}{0.49\textwidth}
\centering
\includegraphics[width=1\linewidth]{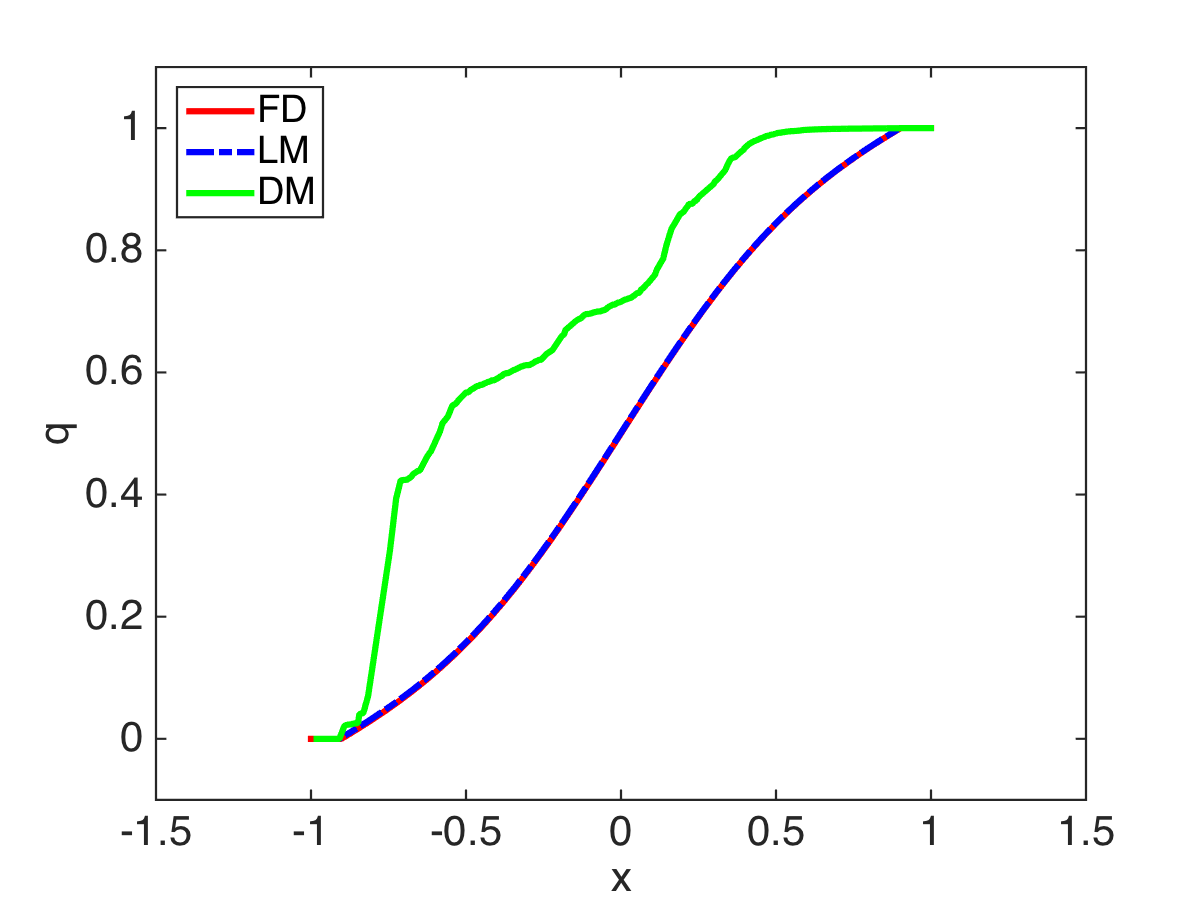}
\end{minipage}
\caption{Left: histogram of distribution of 596 points sampled on $[-1,~1]$. Right: Numerical solutions using different methods.}
\label{fig:Comparison_randpt}
\end{figure}

In the second experiment, we consider points simulated by the SDE
\eqref{eq:overdamp} using an Euler-Maruyama scheme. We take 1000 time snapshots from a long trajectory generated by the SDE~\eqref{eq:overdamp} and only keep those points inside $[-1,~1]$, which provides 596 points. We apply both
local mesh method and diffusion map method to this data set, and
compare to the reference solution obtained by the finite element
method based on $1000$ equally distributed points on $[-1,~1]$. As we
can see from the left panel in Figure \ref{fig:Comparison_randpt},
although the data does not well sample the invariant measure
associated with the double well as we only choose $1000$ points to
discretize the SDE (in particular, the empirical distribution is far
from symmetric), our method still provides almost identical
result as the solution obtained from the standard finite element
method. The diffusion map method on the other hand does not provide
satisfactory result in this case. Note that the information of the invariant measure has been also used in the diffusion map method as the way of constructing $K(\bm{p_i},\bm{p_j})$ in \eqref{eqn:diffKernel}.  

\begin{figure}[h]
\begin{minipage}{0.49\textwidth}
\centering
\includegraphics[width=1\linewidth]{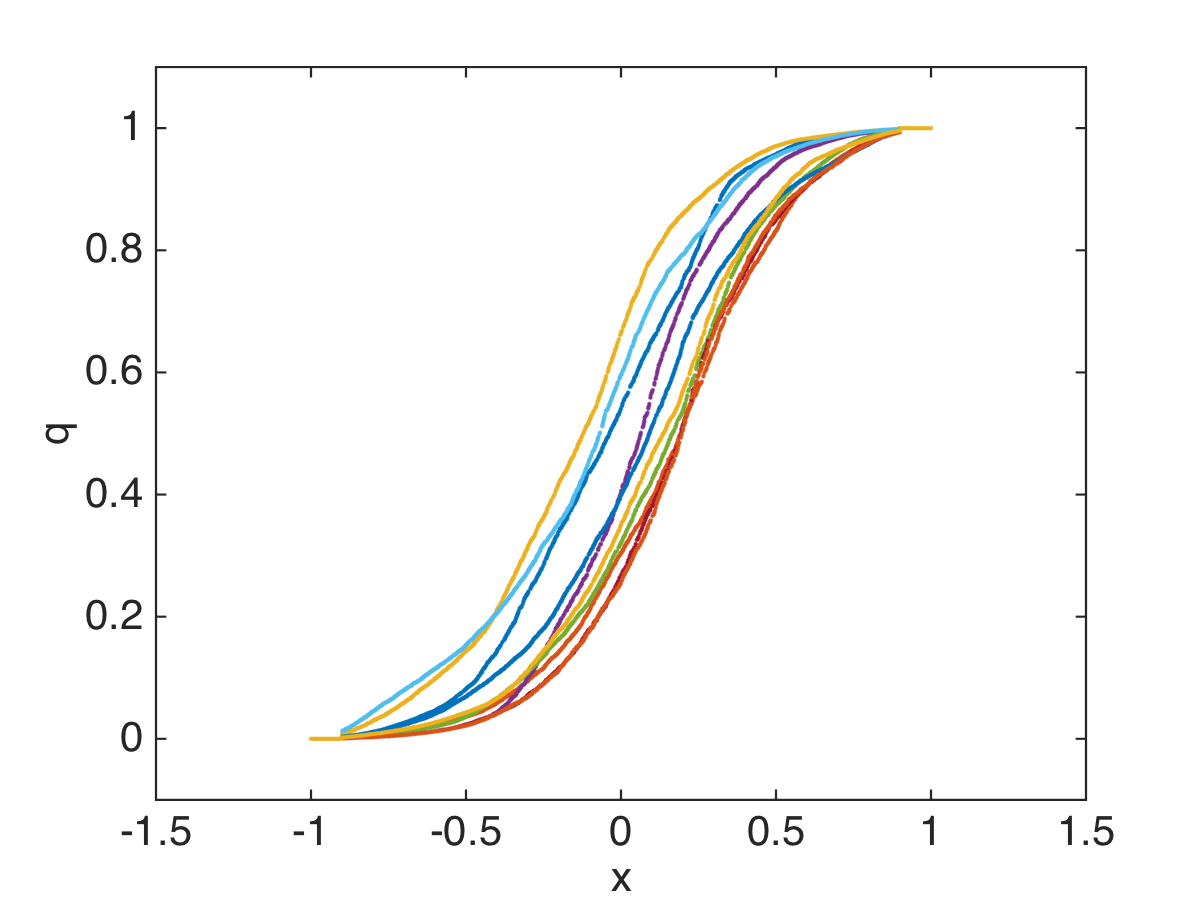}
\end{minipage}
\begin{minipage}{0.49\textwidth}
\centering
\includegraphics[width=1\linewidth]{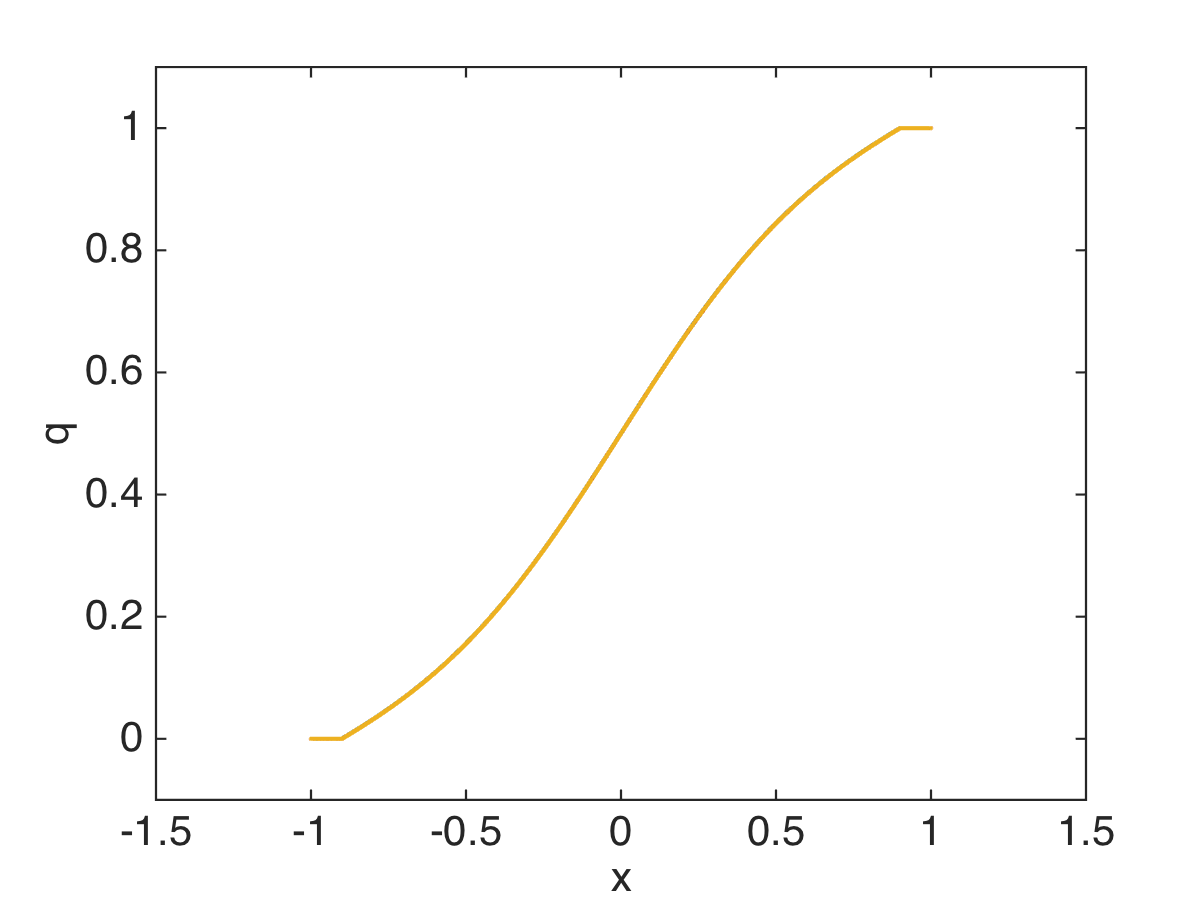}
\end{minipage}
\caption{Left:10 tests to the backward Fokker-Planck equation (sampled from 10000 points) using diffusion map. Right: 10 tests to the backward Fokker-Planck equation (sampled from 10000 points) using our method.}
\label{fig:FPCurve_reprod}
\end{figure}

It is natural to ask about the performance of the methods when more
data points are available; in particular, the accuracy of the
diffusion map discretization will improve. In the third experiment, we
compare our method with diffusion map method using more points
obtained from the same SDE \eqref{eq:overdamp} than the previous
experiment: For each test, we take $10000$ time snapshots from a long
trajectory generated by the SDE and only keep those points inside the
interval $[-1,~1]$, which provides around $6000$ points for each
test. The test is repeated for $10$ times with independent drawing of the point clouds.  Figure
\ref{fig:FPCurve_reprod} reports the approximation of committor
functions as solution of equation \eqref{eq:committor}. The left
panel of Figure \ref{fig:FPCurve_reprod} shows results obtained from
the diffusion map method, where curves with different color-coding
indicate solutions for different realization of the test (note that
the point cloud changes from test to test). It is clear to see that
the diffusion map based method is quite sensitive to sampling quality
of the point clouds, thus the approximated solutions of the committor
function is not quite reproducible due to the different sampling
quality. As a comparison, the corresponding solutions using our method
are plotted in the right picture of Figure \ref{fig:FPCurve_reprod}
using exactly the same samples. We found that solutions
obtained from our method are stacked on top of each other as they are
all essentially identical to the result obtained from the finite
element method based regular grid. This shows
that our method is very robust to points sample quality, and thus
yields better reproducibility, besides provides more accurate
approximation of the committor function.


\subsection{Rugged Mueller potential}
Next we test our methods on a $2D$ problem with the rugged Mueller
potential, which is a well-established test problem in chemical
physics. The point cloud is sampled in $2D$ domain
$\Omega = [-1.5,~1]\times [-0.5,~2]$ and the potential is given by 
\begin{equation}
U_{rm}(x,y) =  U(x,y) + \gamma\sin(2 k \pi x)\sin(2k\pi y)
\label{eqn:ruggedMueller}
\end{equation}
with $U$ being the original Mueller potential 
\begin{equation*}
U(x,y)= \sum_{i=1}^4 D_i e^{a_i (x - X_i)^2 + b_i(x - X_i) (y - Y_i) + c_i (y - Y_i)^2}
\end{equation*}
with the parameters chosen as: 
\begin{align*}
[a_1,a_2,a_3,a_4] & = [-1,-1,-6.5,0.7], \\
[b_1,b_2,b_3,b_4] & = [0,0,11,0.6], \\
[c_1,c_2,c_3,c_4] & = [-10,-10,-6.5,0.7],  \\
[D_1,D_2,D_3,D_4] & = [-200,-100, -170,15], \\
[X_1,X_2,X_3,X_4] & = [1,0 ,-0.5, -1], \\ 
[Y_1,Y_2,Y_3,Y_4] & = [0,~  0.5,~  1.5,~  1].
\end{align*} 
In addition, we choose $\gamma = 9, k = 5$. 
Thus the rugged Mueller potential increases the roughness of the potential.  

\begin{figure}[h]
\centering
\begin{minipage}{0.495\linewidth}
\includegraphics[width=1\linewidth]{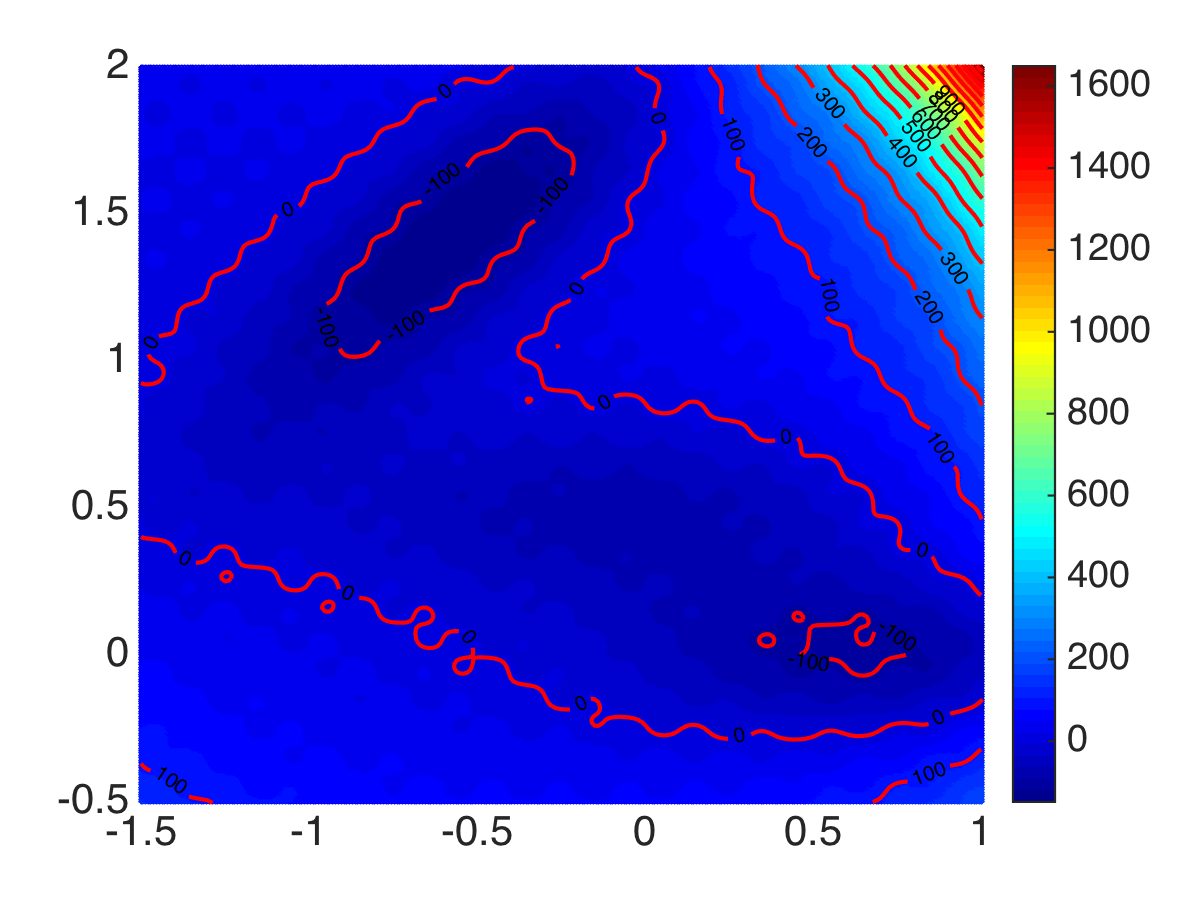}
\end{minipage}
\begin{minipage}{0.495\linewidth}
\includegraphics[width=1\linewidth]{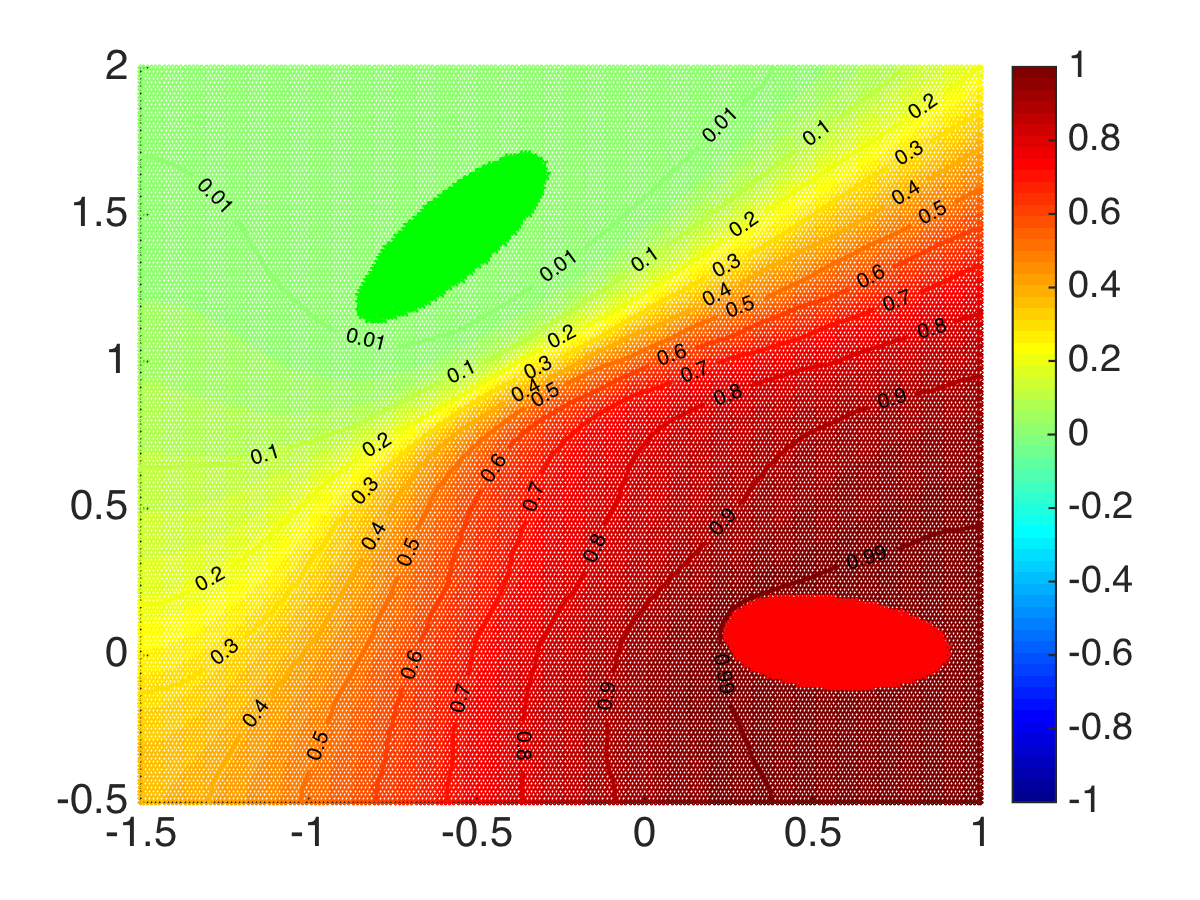}
\end{minipage}
\caption{Left: A rugged Mueller potential used in our
  experiments. Right: The committor function with its level contours
  obtained by the standard finite element method on the domain
  $[-1.5,~1]\times [-0.5,~2]$ , $\nu_R^{\text{FE}} = 0.92960$. }
\label{fig:ruggedMueller}
\end{figure}

The reactant and product sets are chosen to be 
\begin{align*}
  & A = \{ U(x,y) < -120 \} \cap \{ y > 0.75\} \\
  & B = \{ U(x,y) < -82\} \cap \{ y < 0.35 \}.
\end{align*}
In the left
panel of Figure \ref{fig:ruggedMueller}, we plot the above rugged
Mueller potential and its level contours. As a reference, we also use the weak formula \eqref{eqn:weakFK} to
solve the Fokker-Planck equation using the finite element method in the domain
$[-1.5,~1]\times [-0.5,~2]$ and denote the resulting committor function as $q_{\text{FE}}$. The right panel of
Figure~\ref{fig:ruggedMueller} illustrates this committor function and
its level contours. Using the stiffness matrix and the mass matrix constructed in the finite element method, we also numerically evaluate the transition rate $\nu_R$ given by equation~\eqref{eq:nuR}, measures the total probability of reactive trajectories out of
the set $A$ to the set $B$. In this case, we obtain the transition
rate is $0.92960$, which we denote it as $\nu_R^{\text{FE}}$ for later
comparisons. 

\begin{figure}[h]
\begin{minipage}{0.49\linewidth}
\includegraphics[width=1\linewidth]{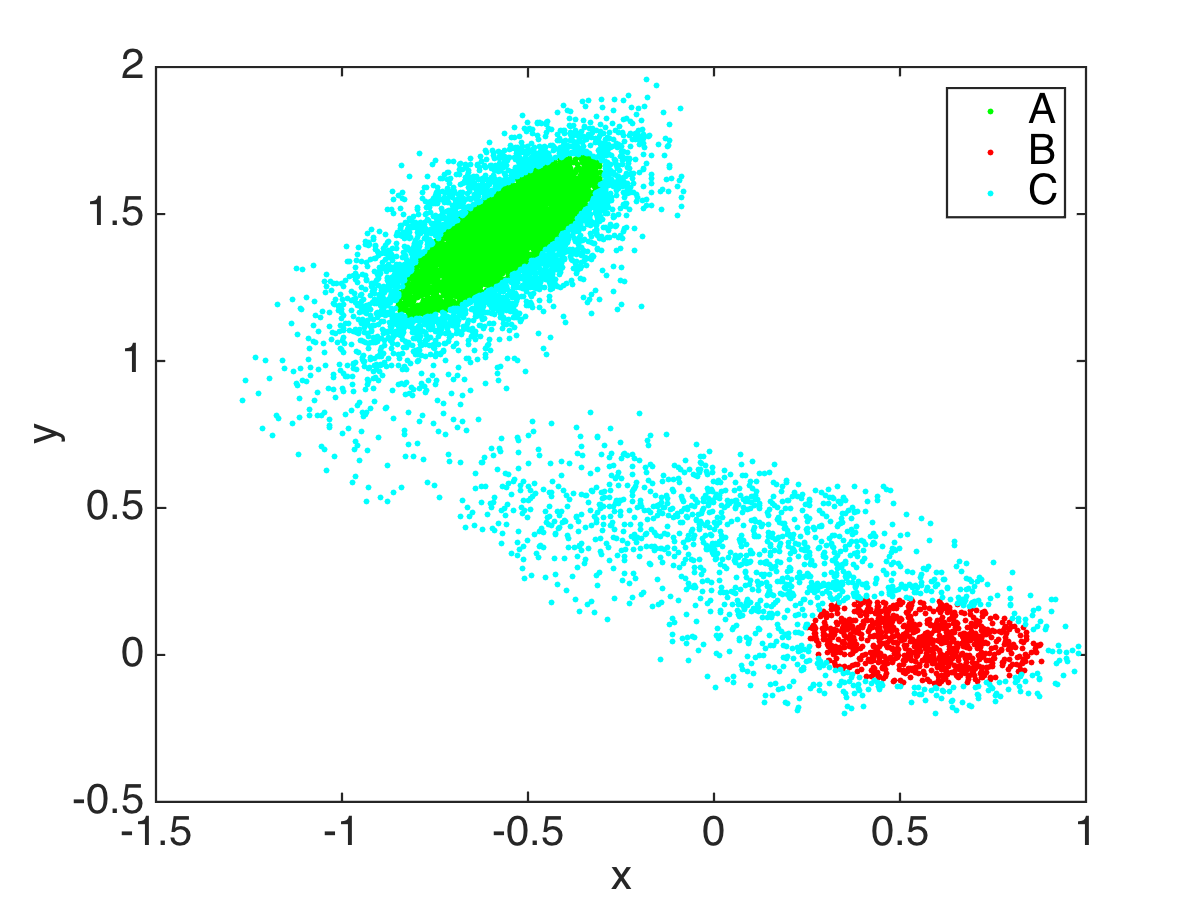}\\
\end{minipage}\hfill
\begin{minipage}{0.49\linewidth}
\includegraphics[width=1\linewidth]{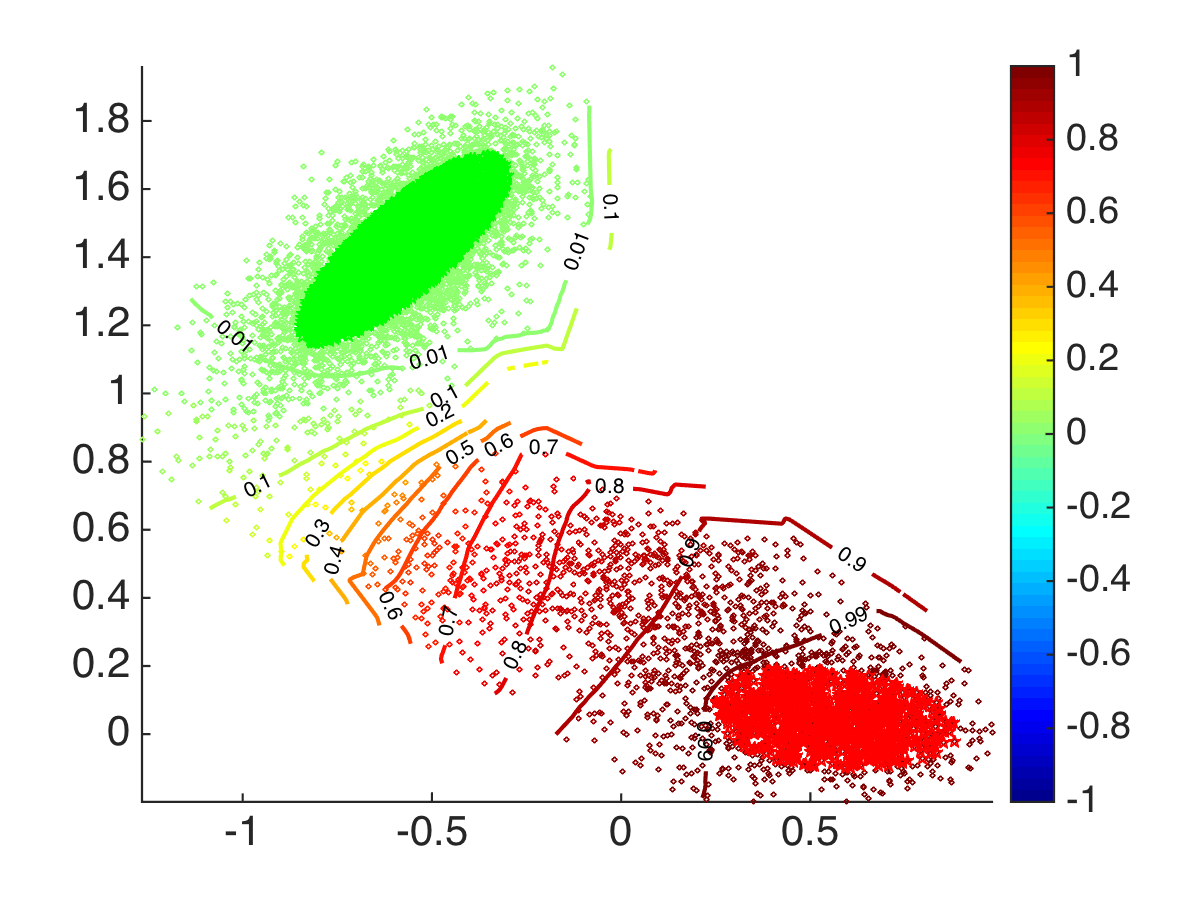}
\end{minipage}\hfill
\caption{Left: Data generated by the SDE \eqref{eq:overdamp} with $\beta = 1/22$. Right: A solution obtained by LM and its level contours.}
\label{fig:2DFPsolution}
\end{figure}

We now consider point clouds generated by numerically integrating the
overdamped Langevin equation \eqref{eq:overdamp} with the above rugged
Mueller potential using the Euler-Maruyama method.  We only keep those
points inside the domain $ [-1.5,~1]\times [-0.5,~2]$. Based on an
input set of irregular data points, the proposed local mesh method
will be applied to solve for the committor function $q$.  In
Figure~\ref{fig:2DFPsolution}, the left panel shows a realization of
the point cloud obtained 
with $\beta = 1/22$. The numerical approximation of a committor
function $q$ is color-coded on the point cloud and illustrated in the
right panel of Figure \ref{fig:2DFPsolution}. Qualitatively, the
obtained function $q$ represents an increase trend of probability that
moving from the set $A$ to the set $B$, consistent with the intuition
behind the committor functions.

To test the accuracy and convergence of the proposed local mesh
method, we chose the number of snapshots in SDE \eqref{eq:overdamp}
 from $20,000$ to $200,000$ with an incremental size
$10,000$. Thus, we obtain a sequence of point clouds which are sampled
in $ [-1.5,~1]\times [-0.5,~2]$ and accumulate near the minimal energy
path connecting two deep wells of the ragged Mueller potential. 

\begin{figure}[h]
\begin{minipage}{0.49\textwidth}
\centering
\includegraphics[width=1\linewidth]{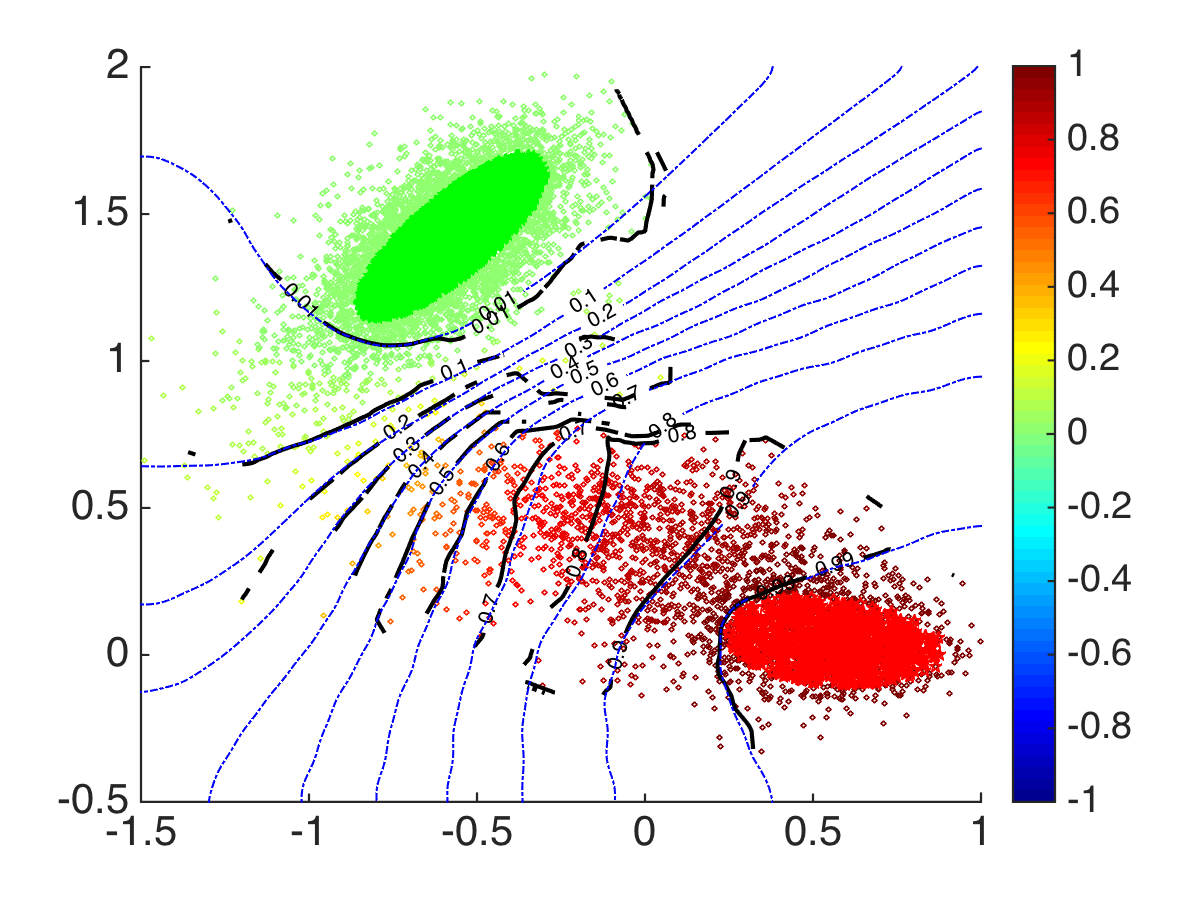}
\end{minipage}\hfill
\begin{minipage}{0.49\textwidth}
\centering
\includegraphics[width=1\linewidth]{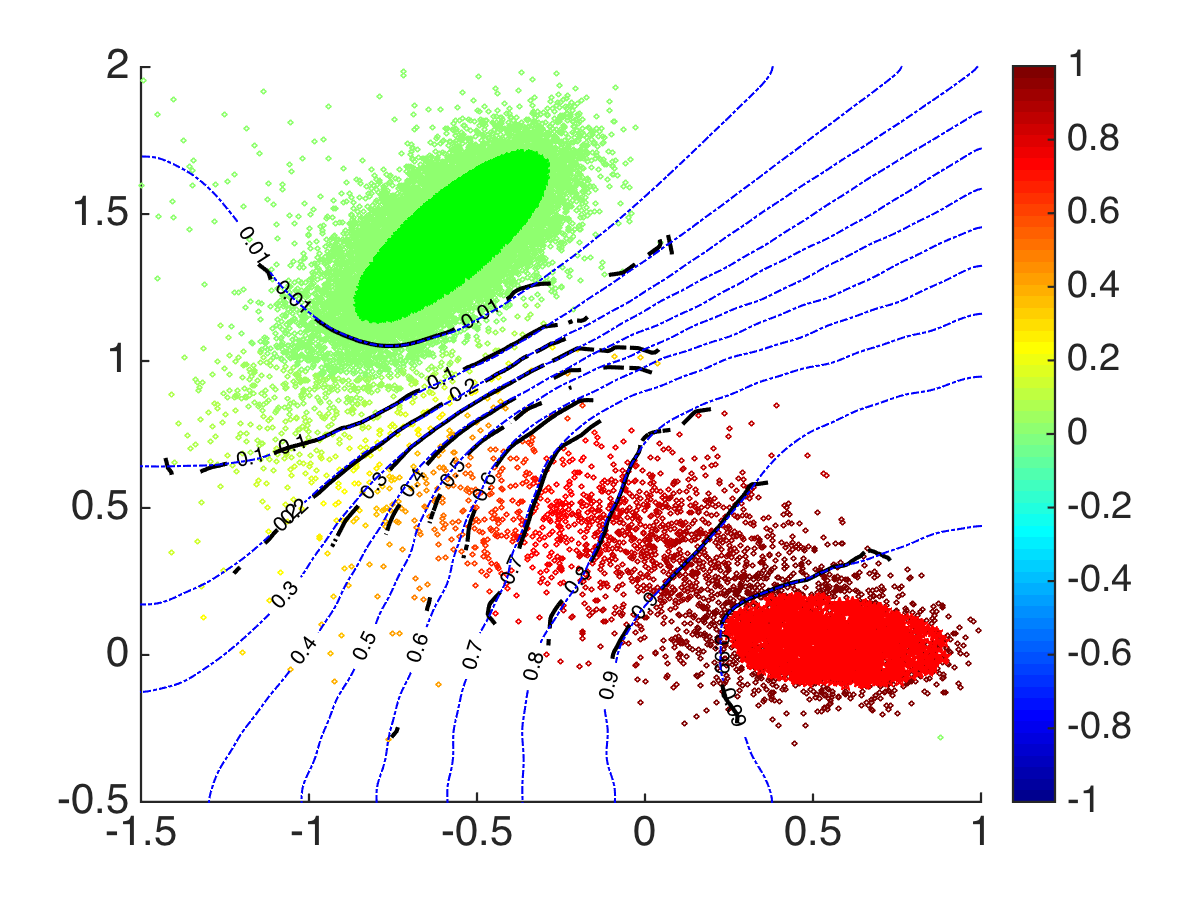}
\end{minipage}\hfill\\
\begin{minipage}{0.49\textwidth}
\centering
\includegraphics[width=1\linewidth]{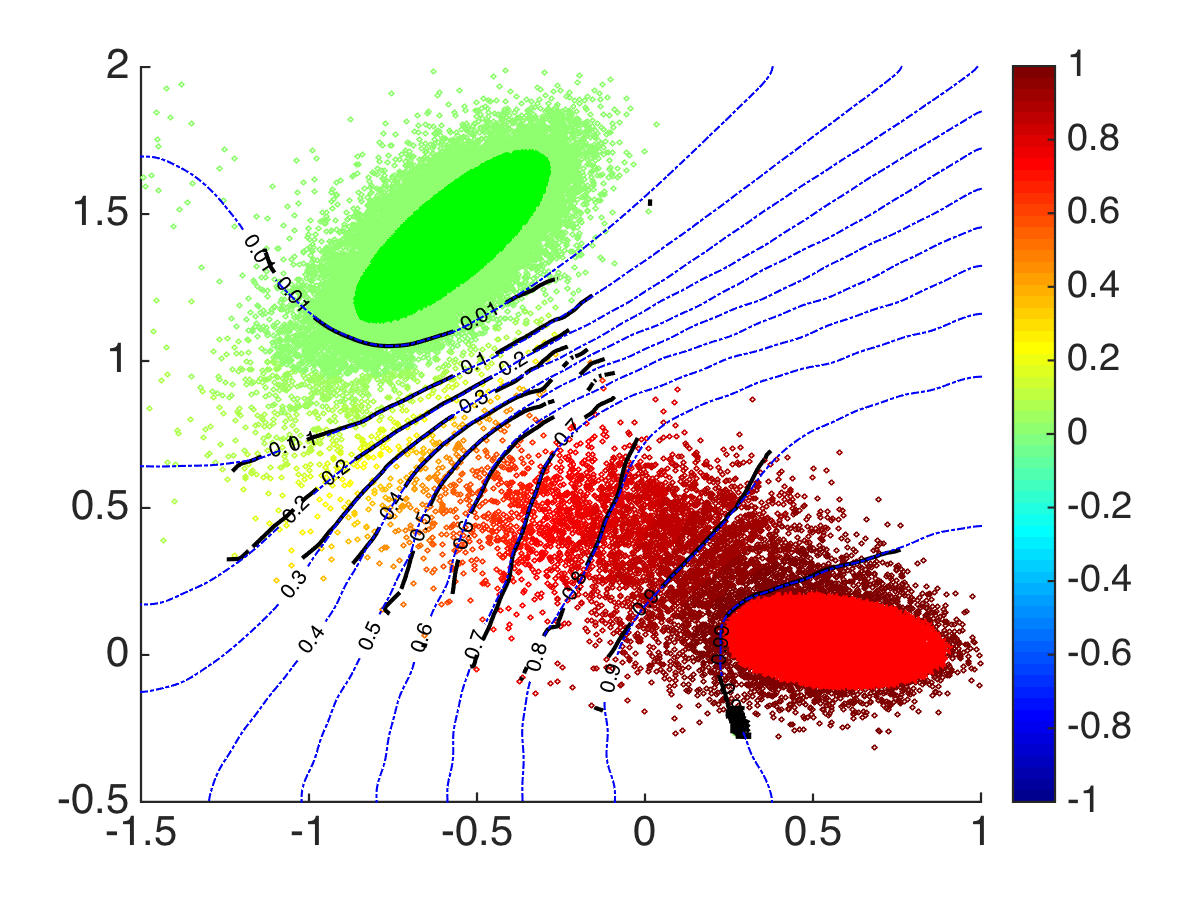}
\end{minipage}\hfill
\begin{minipage}{0.49\textwidth}
\centering
\includegraphics[width=1\linewidth]{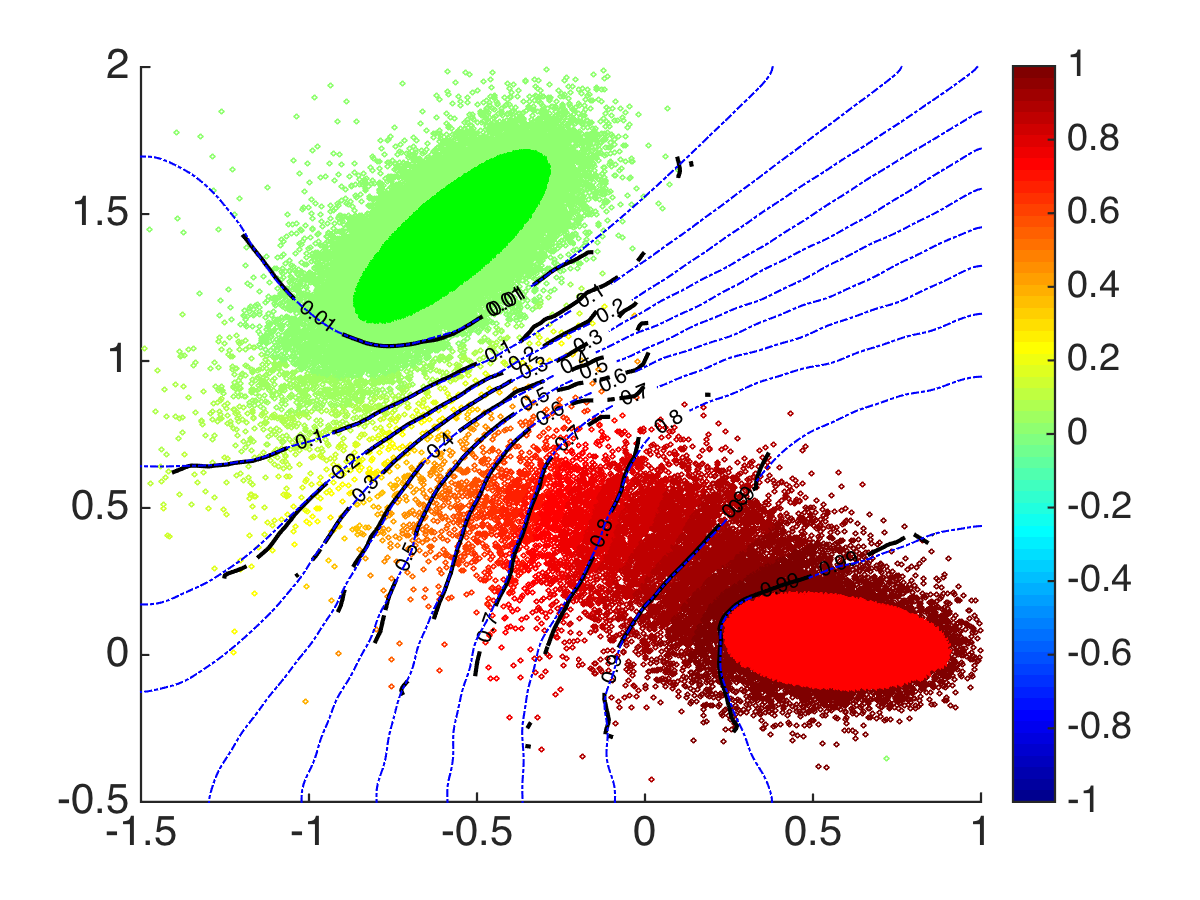}
\end{minipage}
\caption{Contour lines match, where black curves represent contours from our solutions with different number of points generated by SDE \eqref{eq:overdamp} and blue-dash curves represent contours from solutions based on regular grid. Left up: 11372 points. Right up: 25711 points. Left down: 45035 points. Right down: 69274 points.}
\label{fig:Contours}
\end{figure}

\begin{figure}[h]
\begin{minipage}{0.49\textwidth}
\centering
\includegraphics[width=1\linewidth]{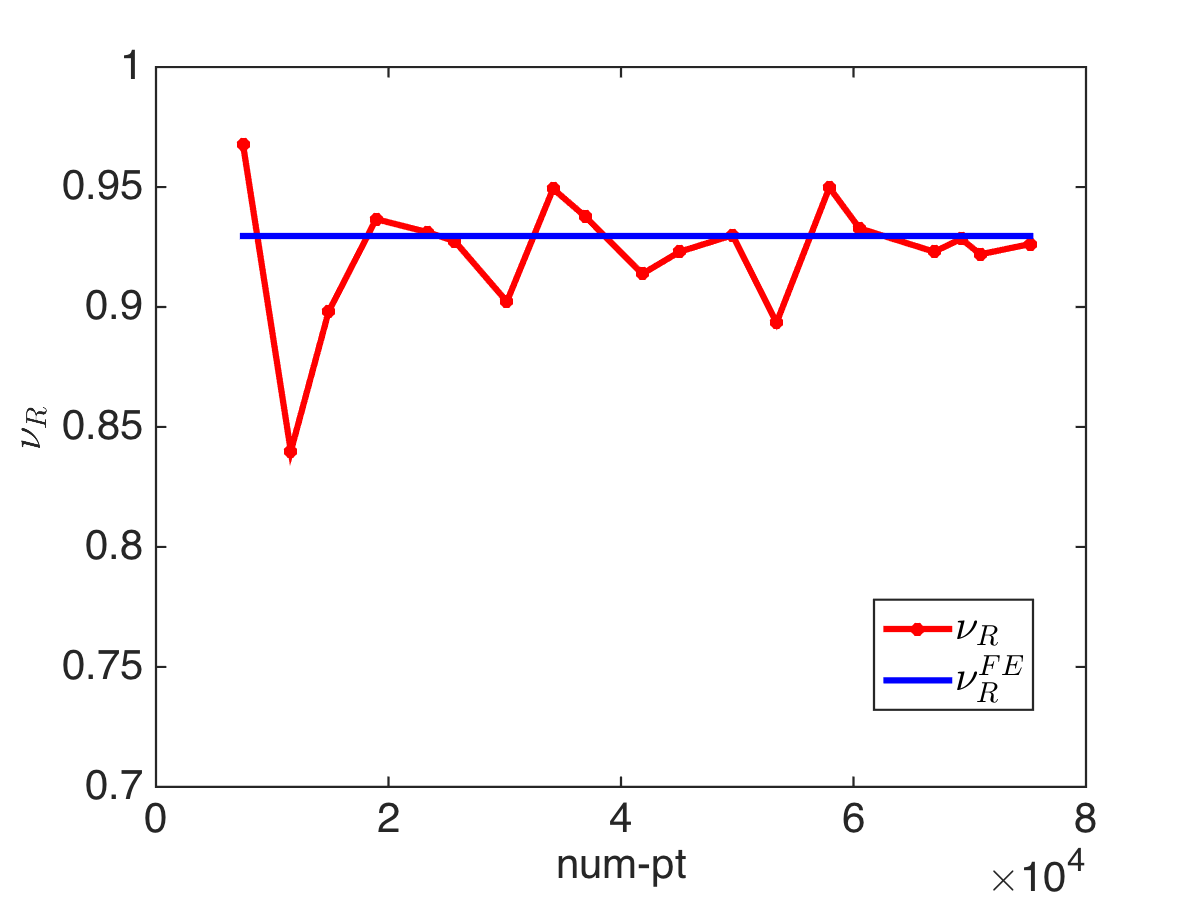}
\end{minipage}
\begin{minipage}{0.49\textwidth}
\centering
\includegraphics[width=1\linewidth]{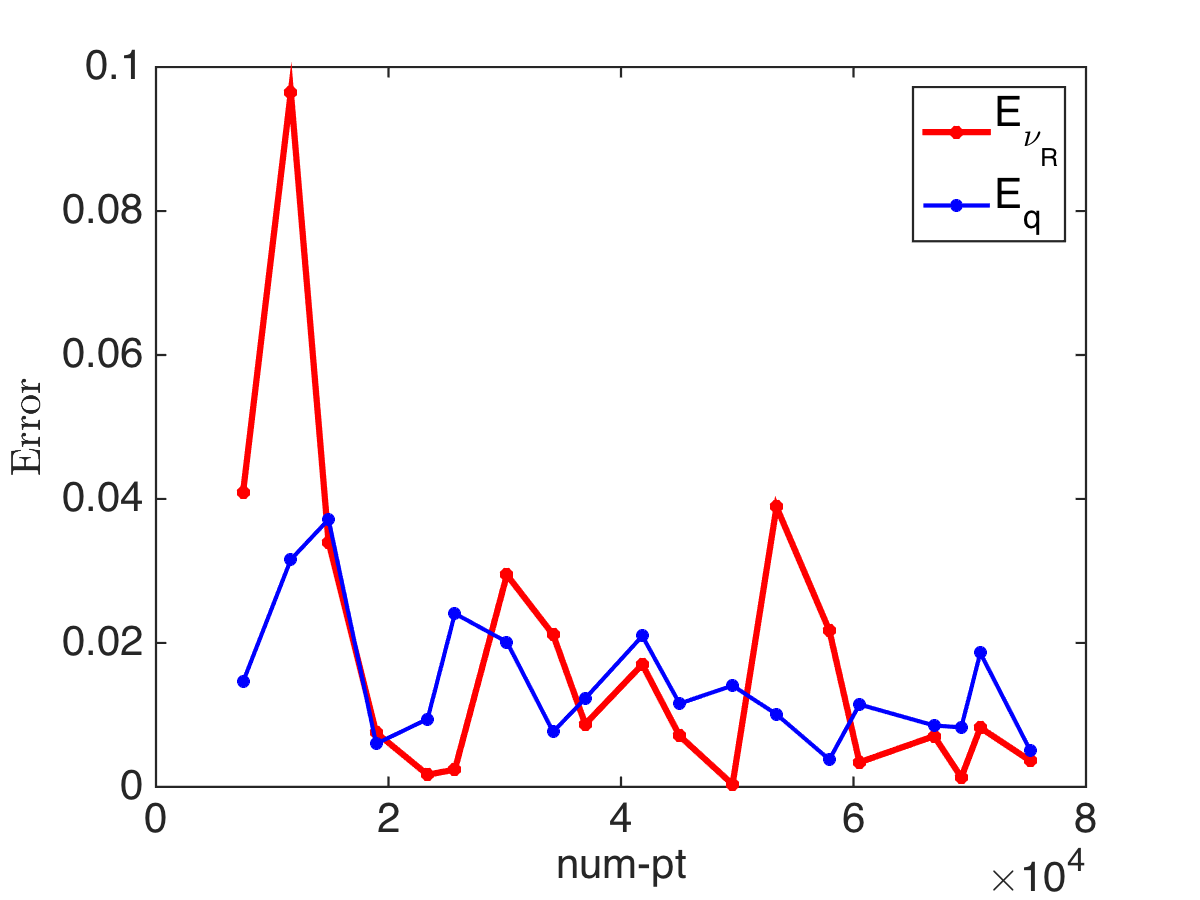}
\end{minipage}
\caption{Left: $\nu_R$ v.s number of points. Right: Relative error of $\nu_R$ and relative error of $q$ v.s. number of points. }
\label{fig:VrConvergence}
\end{figure}

The local mesh method thus produces a sequence approximation of the committor function. For
comparison, we use committor function $q_{\text{FE}}$ based on the regular grid
$[-1.5,~1]\times [-0.5,~2]$ as a reference. As can be seen from Figure~\ref{fig:Contours}, level contours, represented by black curve, of the
approximated committor function match very well to the blue curves
representing level contours of $q_{\text{FE}}$. Moreover, the transition rate $\nu_R$ defined in ~\eqref{eq:nuR} can be approximated by evaluating $\int_{\Omega} |\nabla q| e^{\beta U} \mathrm{d}x$ using the stiffness matrix constructed in~\eqref{eqn:Stiffness2} and approximating $Z = \int_{\Omega} e^{-\beta U} \mathrm{d}x$ using the method proposed in~\cite{LaiLiangZhao:13} based on a numerical approximation of the mass matrix. After that, we measure the
relative error
\[ 
E_{\nu_R} = \dfrac{|\nu_R -\nu_R^{\text{FE}}|}{\nu_R^{\text{FE}}}
\]
between the approximated transition rate $\nu_R$ and the transition
rate $\nu_R^{\text{FE}}$. We also compute the relative error
\[
E_{q} = \frac{\|q - \tilde{q}_{\text{FE}}\|_2}{\|\tilde{q}_{\text{FE}}\|_2}
\]
 between $q$ and $q_{\text{FE}}$. We remark that the original $q_{\text{FE}}$ is defined on regular grid. 
 The above $\tilde{q}_{\text{FE}}$ is calculated using the interpolation
of $q_{\text{FE}}$ from regular grid to the input points and the standard $l_2$ norm is used here to measure difference between these two vectors. Our numerical results reported in Figure
\ref{fig:VrConvergence} indicate that the approximation error can be
controlled around $1\%$ for moderate size of points.

As a direct application based on the committor function obtained from
the Fokker-Planck equation, we trace a deterministic reactive flow
$X(s)$ by solving the following ODE on the point clouds
\begin{equation}
\label{eqn:trajectorytracking}
\left\{\begin{array}{l}
\displaystyle \frac{dX(s)}{ds} = J_R(X(s)) = k_B T \rho(X(s)) \nabla q(X(s)),\vspace{0.2cm} \\
X(0)= p_0
\end{array}\right.
\end{equation}
The idea of solving the above ODE is to interpolate the point cloud locally using the moving least square method~\cite{Liang:CVPR2012,Liang2013solving}, then the solution curve can be extended based on the locally interpolated manifold. 

\begin{wrapfigure}{r}{0.45\textwidth}
\vspace{-1cm}
\begin{center}
\includegraphics[width=.9\linewidth]{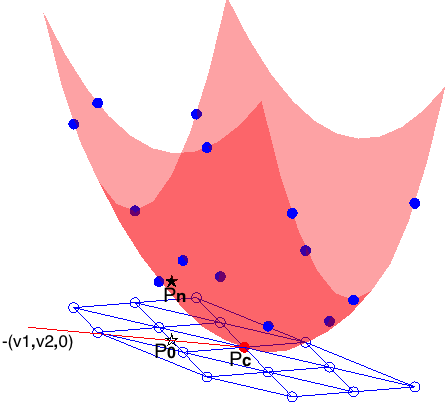} \\
\end{center}
\vspace{-.7cm}
\caption{\small Reactive flow tracing. The current point $\bm{p}_c$, its KNN points and the new point $\bm{p}_{new}$ are marked as the red solid circle, the blue solid circles and the black solid star respectively, whose projections on the tangle plane at $\bm{p}_c$ are plotted as the corresponding hollow markers. The red line has direction $(-v_1,-v_2,0)$.}
\label{fig:ReactiveFlow}
\end{wrapfigure}
\begin{figure}[h]
\centering
\includegraphics[width=.7\linewidth]{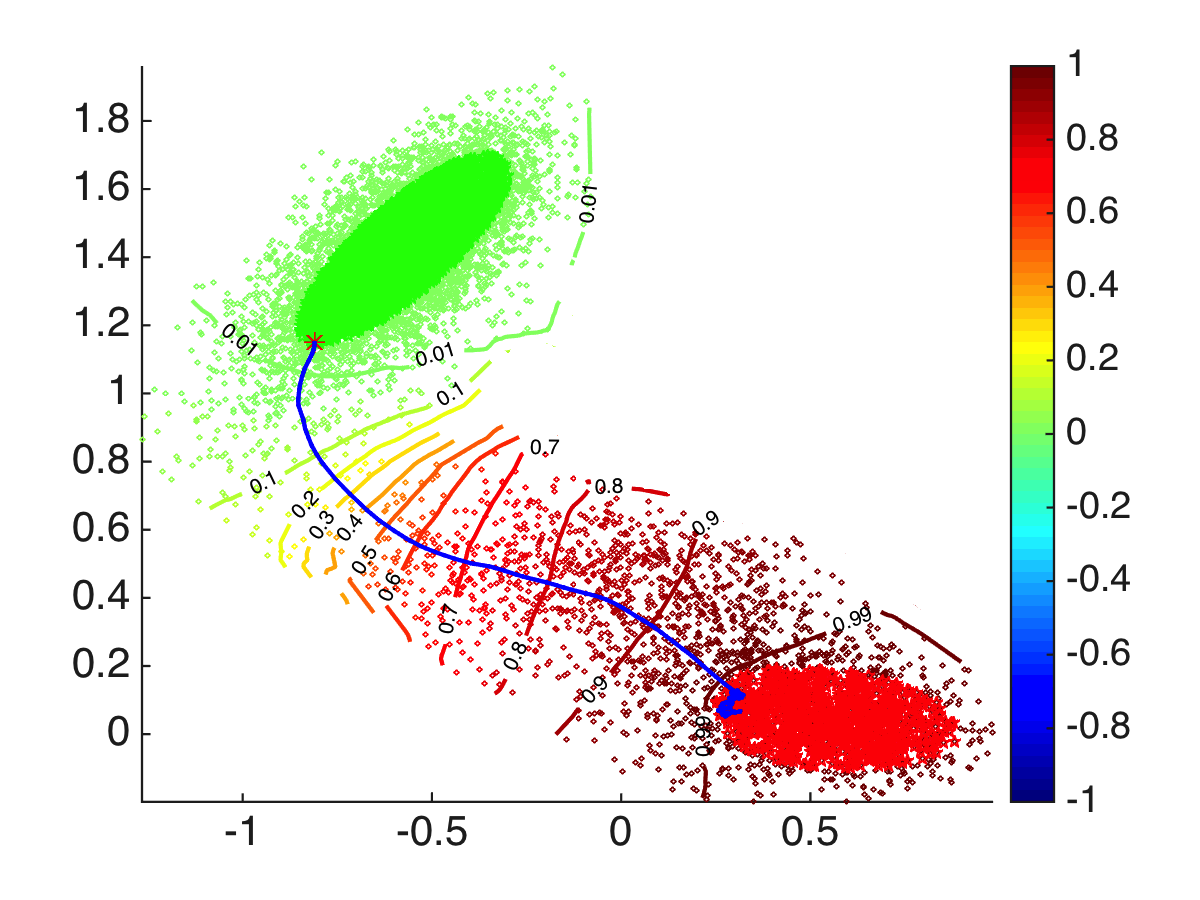}
\caption{Reactive trajectory tracing.}
\label{fig:trajectory}
\end{figure}

To clearly indicate our method of solving \eqref{eqn:trajectorytracking}, we assume that the given point cloud is sampled on a two dimensional manifold embedded in $\RR^3$. We emphasize that the following idea can be straightforwardly extended to high dimension cases. Suppose a point $\bm{p}_c$ (current point on the reactive flow) has already been obtained, we intend to find the next point on the reactive flow. Without loss of generality, suppose $\bm{p}_1,\bm{p}_2, \cdots, \bm{p}_K \in\P$ are KNN of $\bm{p}_c$ in the point cloud $\P$. Using PCA, we can build a local coordinate system $\{\bm{e}_{\bm{p}_c}^1, \bm{e}_{\bm{p}_c}^2, \bm{e}_{\bm{p}_c}^3\}$ centered at $\bm{p}_c$ and the KNN of $\bm{p}_c$ has local coordinates $(x_i,y_i,z_i)$. We use moving least squares (MLS) to locally approximate the surface as $\Gamma=(x,y,z(x,y))$ and estimate $k_B T \rho(\bm{p}_c)\ \nabla q(\bm{p}_c)$ (more details can be found in \cite{Liang:CVPR2012,Liang2013solving}). We construct the Delaunay triangulation of the projections $\{\hat{\bm{p}}_c,\hat{\bm{p}}_1,\hat{\bm{p}}_2,\cdots,\hat{\bm{p}}_K\}$ and find the first ring $\R=\{T_c^1,\cdots, T_c^l\}$ of $p_c$, which is the same as we did in Section \ref{subsec:localmesh}. Suppose that $k_B T \rho(\bm{p}_c)\ \nabla q(\bm{p}_c)$ has a local coordinate $(v_1,v_2,v_3)$ in $\{\bm{e}_{\bm{p}_c}^1,\bm{e}_{\bm{p}_c}^2,\bm{e}_{\bm{p}_c}^3 \}$, 
we find the intersection of line segment starting at $\bm{p}_c$ with the direction $(-v_1,-v_2,0)$ and the first ring. Notice that this computation is done within the tangent space of $\bm{p}_c$. Denote the intersection as $\hat{\bm{p}}_{new}=(x_0,y_0,0)$, we then project it back to the approximated surface to obtain the next point on the geodesic path $\bm{p}_{new}=(x_0,y_0,z(x_0,y_0))$. This process is illustrated in Figure \ref{fig:ReactiveFlow}. We refer~\cite{LaiLiangZhao:13} for more detailed discussion about solving the above equation on point clouds. Figure~\ref{fig:trajectory} plots the trajectory starting from the red star point $p_0$ in state $A$ to finally hit the region in state $B$, which clearly show that the reactive flow jump from state $A$ to state $B$. 

\subsection{Experiments in higher dimensions}
\begin{figure}[h]
\begin{minipage}{0.49\textwidth}
\centering
\includegraphics[width=1\linewidth]{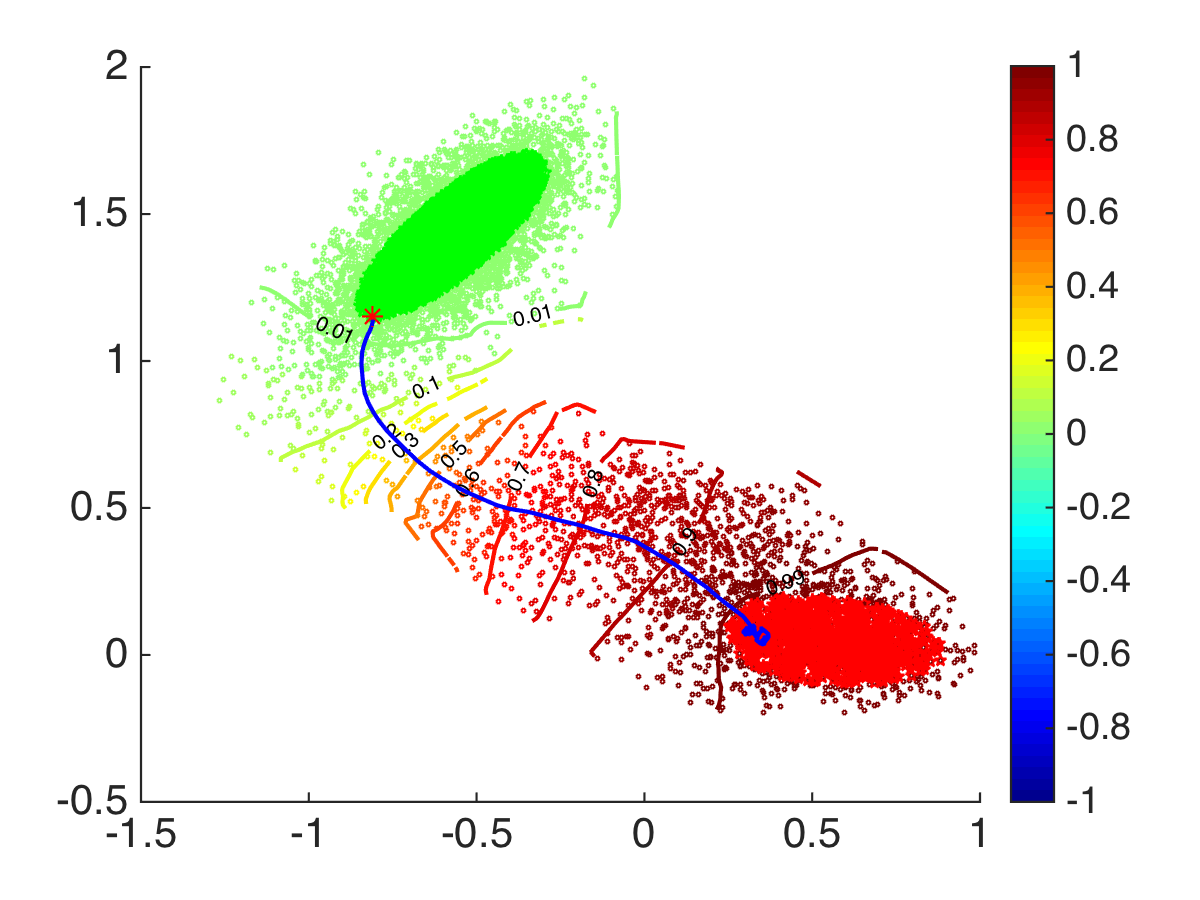}
\end{minipage}\hfill
\begin{minipage}{0.49\textwidth}
\centering
\includegraphics[width=1\linewidth]{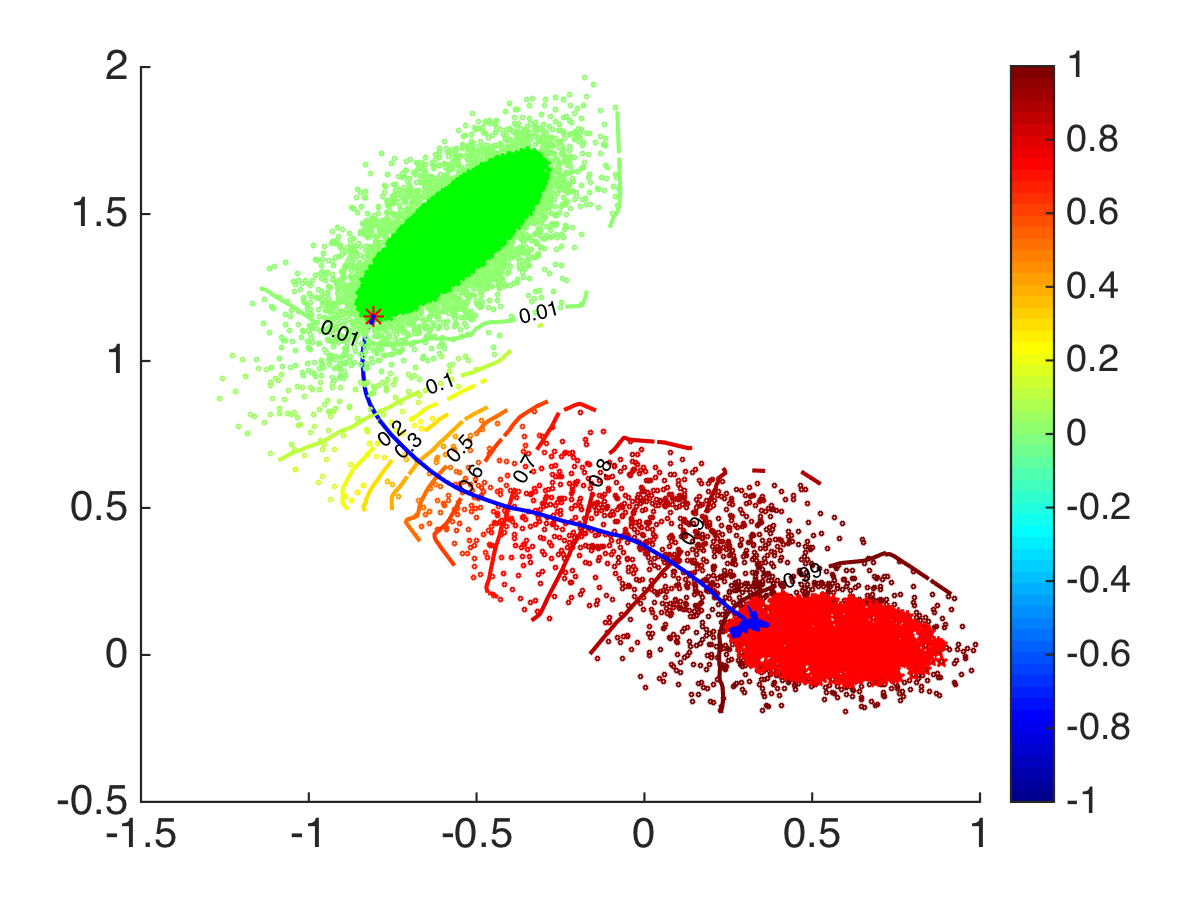}
\end{minipage}\hfill\\
\begin{minipage}{0.49\textwidth}
\centering
\includegraphics[width=1\linewidth]{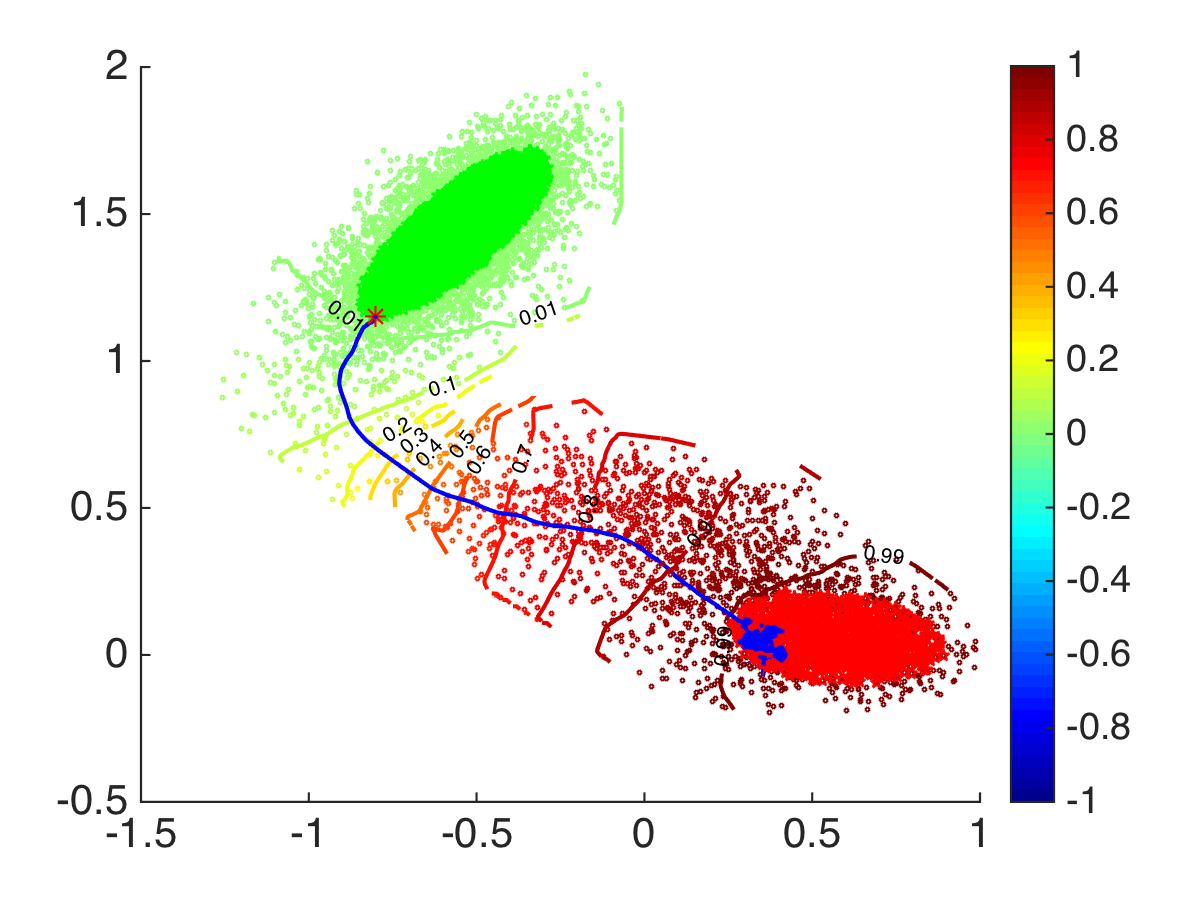}
\end{minipage}\hfill
\begin{minipage}{0.49\textwidth}
\centering
\includegraphics[width=1\linewidth]{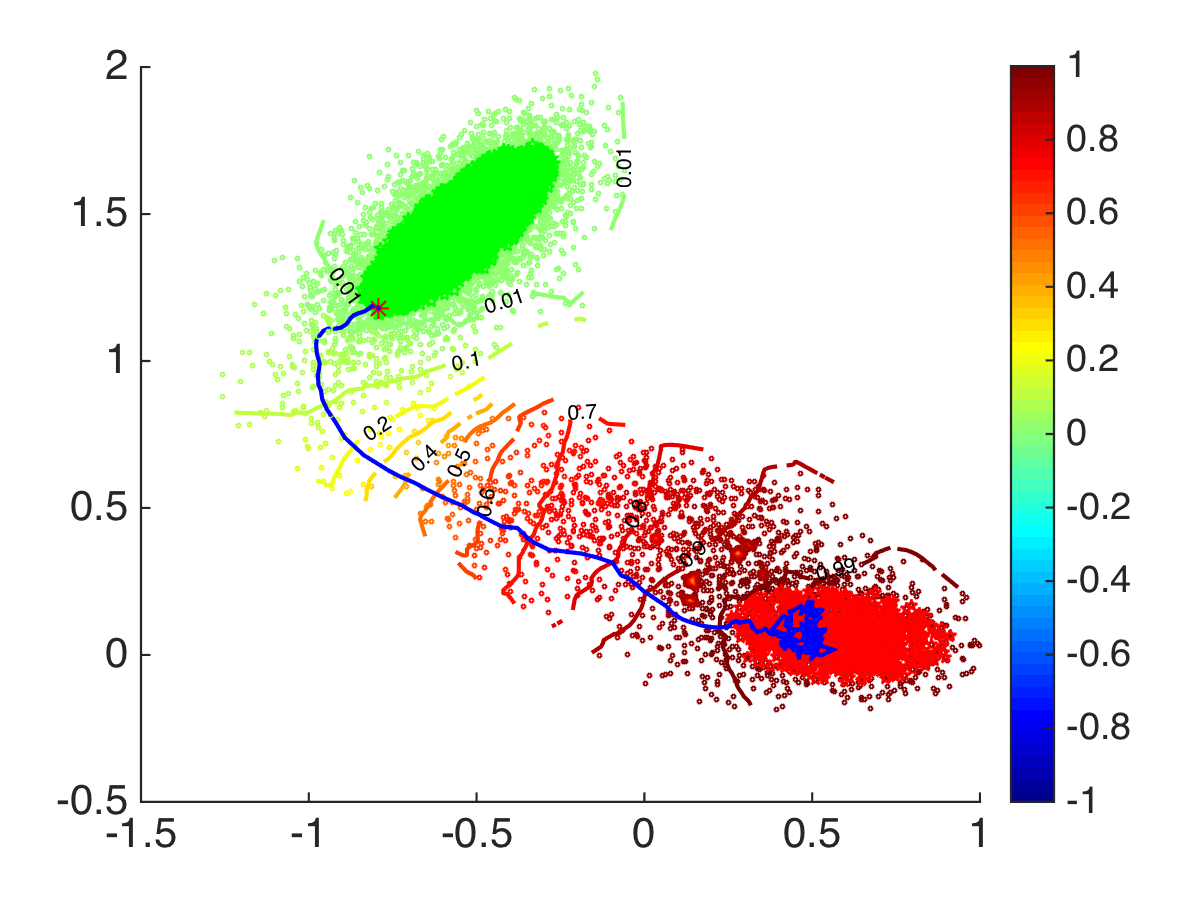}
\end{minipage}
\caption{The committor functions and reaction trajectories for point clouds embedded in $\RR^{10}$ with noise levels $\gamma = 10\%, 20\%, 50\%, 100\%$ from top left to bottom right, respectively.}
\label{fig:NoiseHD}
\end{figure}

As an advantage of the proposed intrinsic method, the solver we designed can handle point clouds sampled from a low-dimensional manifold in a high dimension space. To illustrate the robustness of the proposed method, we next test the our solver for point clouds embedded in a high dimensional ambient space with an artificial Gaussian noise. In other words, we first simulate 2D point clouds as we conducted in the previous numerical experiments. After that, we embed the point cloud to $\RR^{10}$ by setting the last eight coordinates to be zero. In addition, we also perturb this point cloud in $\RR^{10}$ by adding Gaussian noise with variance $\sigma = \gamma d_{max}$. Here we choose $d$ to be the maximal number of the 50th smallest distance to each point. Namely, $d_{max} = \max_i\{d_i~|~ d_i = \mbox{the 50-th smallest distance to } \bm{p}_i\}$. In our experiments, we solve the Fokker-Planck equation based on the point cloud perturbed by Gaussian noise with different level $\gamma = 10\%, 20\%, 50\%, 100\%$, and also compute the reactive trajectory from the starting point. For better visualization, figure~\ref{fig:NoiseHD} shows 2D projection of the Gaussian noise perturbed point clouds color-coded with the resulting commitor function $q$. In addition, reactive trajectories are also plotted for different noise level with the same starting point. This figure clearly demonstrates the robustness of the proposed method. 

\section{Conclusions}
\label{sec:conclusions}

In this work, we develop a point cloud discretization for computing
committor functions of stochastic systems. Numerical examples on toy
model systems confirm that the method provides a promising tool to
analyze the stochastic system in the framework of the transition path
theory. In particular, the point cloud discretization extends the
applicability of the transition path theory beyond the ``tube
approximation''. As for future directions, an obvious next step is to
test the approach in thermally activated process in more complicated
and realistic examples arising from biophysics. In addition, our
method does not require that the point cloud samples exactly the
invariant measure. This provides advantages for considering point
clouds sampled locally rather than using a long trajectory and for
combining our method with advanced sampling strategies of the underlying
stochastic system.  The numerical convergence analysis of the point
cloud discretization for Fokker-Planck operators is also an
interesting topic to pursue.



\section*{Acknowledgments}
This work was supported in part by the National Science Foundation
through the grants DMS-1522645 (R.L.) and DMS-1454939
(J.L.). Authors would like to thank Mauro Maggioni for helpful
discussions in various stages of this project.


\end{document}